\documentclass[sn-mathphys]{sn-jnl-mod}

 \normalbaroutside

\theoremstyle{thmstyleone}%
\newtheorem{thm}{Theorem}
\newtheorem{prop}{Property}
\newtheorem{problem}{Problem}
\newtheorem{rem}{Remark}

\newcommand{\bo}[1]{{\bf #1}}
\graphicspath{{./pictures/}}

\title{Parametric Shape Optimization using the Support Function}

\author[1]{\fnm{Pedro R.S.} \sur{Antunes}}\email{prantunes@fc.ul.pt}

\author[2]{\fnm{Beniamin} \sur{Bogosel}}\email{beniamin.bogosel@polytechnique.edu}

\affil[1]{\orgdiv{Sec\c{c}\~{a}o de Matem\'{a}tica, Departmento de Ci\^{e}ncias e Tecnologia}, \orgname{Universidade Aberta}, \orgaddress{\street{Rua da Escola Polit\'{e}cnica 141-147}, \city{Lisbon}, \postcode{1269-001}, \country{Portugal}}}

\affil[2]{\orgdiv{Centre de Math\'ematiques Appliqu\'ees}, \orgname{\'Ecole Polytechnique}, \orgaddress{\street{Rue de Saclay}, \city{Palaiseau}, \postcode{91128}, \country{France}}}

\begin{document}

\title{Parametric Shape Optimization using the Support Function}

\abstract{
The optimization of shape functionals under convexity, diameter or constant width constraints shows numerical challenges. The support function can be used in order to approximate solutions to such problems by finite dimensional optimization problems under various constraints. We propose a numerical framework in dimensions two and three and we present applications from the field of convex geometry. We consider the optimization of functionals depending on the volume, perimeter and Dirichlet Laplace eigenvalues under the aforementioned constraints. In particular we confirm numerically Meissner's conjecture, regarding three dimensional bodies of constant width with minimal volume.
}


\maketitle 

\section{Introduction} 

Shape optimization problems are a particular class of optimization problems where the variable is a shape. A typical example of such a problem has the form
\[
\min_{\omega \in \mathcal{A}} \mathcal{F}(\omega),
\]
where the functional $\mathcal{F}$ is computed in terms of the shape $\omega$ and $\mathcal{A}$ is a family of sets with given properties and eventual constraints. The dependence of the cost functional $\mathcal{F}$ on the geometry can be explicit (volume, perimeter) or implicit, via a partial differential equation. Classical examples in this sense are functionals depending on the spectrum of various operators related to the shape $\omega$, like the Dirichlet-Laplace operator.

When dealing with constrained shape optimization problems, having volume or perimeter constraints facilitates the study of optimizers, in particular because there exist arbitrarily small inner and outer perturbations of the boundary which preserve the constraint. This is not the case anymore when working in the class of convex sets or when bounds on the diameter are imposed. The papers \cite{LNconv09},\cite{LNPconv12} describe some of the theoretical challenges when working with these constraints. 

Challenges of the same nature arise when dealing with convexity, constant width and diameter constraints from a numerical point of view. There are works in the literature which propose algorithms that can handle the convexity constraint. In \cite{LROconvex} a convex hull method is proposed in which the convex shapes are represented as intersections of half-spaces. In  \cite{MOconvex} the authors propose a method of projection onto the class of convex shapes. The articles \cite{BLOcw}, \cite{spheroforms}, \cite{oudetCW} show how to deal with width constraints. The approach proposed in \cite{bartels-wachsmuth} handles simultaneously convexity and PDE constraints by considering discretized domains (typically triangulations) and using deformations which preserve convexity. The methods presented in the previous references are rather complex and not straightforward to implement. In this article we present a more direct approach using the properties of the support function. Such a method was already proposed in \cite{BHcw12} for the study of shapes of constant width, but was essentially limited to the two dimensional case. In particular, the three dimensional computations presented there are for rotationally symmetric shapes, which allows the use of two dimensional techniques. Moreover, the numerical framework in \cite{BHcw12} needs special tools regarding semi-definite programming algorithms and the cost functional is at most linear or quadratic in terms of the Fourier coefficients of the support function.

The precise definition and main properties of the support function are presented in Section \ref{support}. Recall that for a convex body $K \subset \Bbb{R}^d$ the associated support function $p$ is defined on the unit sphere $\Bbb{S}^{d-1}$ and for each $\theta \in \Bbb{S}^{d-1}$, $p(\theta)$ measures the distance from a fixed origin, which can be chosen inside $K$, to the tangent hyperplane to $K$ orthogonal to $\theta$. Already from the definition it can be noted that the quantity $p(\theta)+p(-\theta)$ represents the diameter or width of the body $K$ in the direction parallel to $\theta$. This allows to easily transform diameter or constant width constraints into functional inequality or equality constraints in terms of the support function. Convexity constraints can be expressed in similar ways, with complexity varying in terms of the dimension $d$. These facts are recalled in the following section.

Finite dimensional approximations of convex bodies can be built using a truncation of a spectral decomposition of the support function: Fourier series decomposition for $d=2$ and spherical harmonic decomposition for $d = 3$. The shape optimization problem becomes a parametric optimization problem in terms of the coefficients of the spectral decomposition. The convexity constraint is modeled by a set of linear pointwise inequalities for $d=2$ or quadratic pointwise inequalities for $d=3$. The constant width constraint is characterized by the fact that coefficients of the even basis functions are zero. Diameter constraints can also be translated into a set of pointwise linear inequalities. It can be noted that in some particular cases, functionals like volume and perimeter have explicit formulas in terms of the coefficients in the above decompositions.

In \cite{BHcw12} the authors study numerically optimization problems under constant width constraint in dimension two, with the aid of the support function and Fourier series decomposition. They work with a global parametrization of the convexity constraint, which requires the use of specific semidefinite-programming techniques and software. We choose to work in a simplified framework, inspired from \cite{antunesMM16}, in which the convexity constraint is imposed on a finite, sufficiently large, number of points distributed on the unit circle ($d=2$) or on the unit sphere ($d=3$). This gives rise to simpler constrained optimization problems that can be handled using standard optimization software.

In Section \ref{sec:theory} we provide existence results for the problems we consider. Moreover, we prove that solutions obtained when using a truncated spectral decomposition of the support function converge to the solutions of the original problems as the number of coefficients goes to $+\infty$.

The paper is organized as follows. In Section \ref{support} various properties of the support function parametrization in dimension two and three are recalled. In Section \ref{sec:theory} theoretical aspects regarding the existence of solutions and the convergence of the discrete solutions are investigated. Section \ref{sec:numerics} deals with the parametric representation of shapes using the spectral decomposition of the support function. The handling of convexity, constant width, diameter and inclusion constraints is discussed. Section \ref{sec:fsol} recalls the method of fundamental solutions used for solving the Dirichlet-Laplace eigenvalue problems.

Section \ref{sec:applic} contains applications of the numerical framework proposed for various problems in convex geometry. In particular, a confirmation of the Meissner conjecture regarding bodies of constant width with minimal volume in dimension three is provided. The two different Meissner bodies are obtained by directly minimizing the volume under constant width constraint, starting from general random initializations. Further  applications presented in Section \ref{sec:applic} concern the minimization of eigenvalues of the Dirichlet-Laplace operator under convexity and constant width constraints, approximation of rotors of minimal volume in dimension three, approximation of Cheeger sets and the minimization of the area under minimal width constraint.

 The goal of this paper is to present a general method for performing shape optimization under various non-standard constraints: convexity, fixed width, diameter bounds by transforming them into algebraic constraints in terms of a spectral decomposition of the support function. In order to illustrate the method, various numerical results are presented, some of which are new and are listed below:
	
\begin{itemize}
	\item \emph{optimization of the Dirichlet-Laplace eigenvalues under convexity constraint}: the case $k=2$ in dimension two was extensively studied (see for example \cite{antunes_henrot}, \cite{oudeteigs}). The case $k \geq 3$ in dimension two and the simulations in dimension three are new and are presented in Section \ref{sec:eigs}.
	\item \emph{minimization of the Dirichlet-Laplace eigenvalues under fixed constant width constraint in dimension three}: in Section \ref{sec:eigs-cw} cases where the ball is not optimal are presented.
	\item \emph{numerical confirmation of Meissner's conjecture}: in Section \ref{sec:meissner} the proposed numerical framework allows to obtain Meissner's bodies when minimizing the volume under fixed width and convexity constraints, starting from randomized initial spherical harmonics coefficients.
	\item \emph{rotors of minimal volume in dimension three}: the numerical framework allows to approximate optimal rotors in the regular tetrahedron and the regular octahedron. (Section \ref{sec:rotors})
\end{itemize}
Other results deal with the computation of Cheeger sets for various two dimensional and three dimensional domains and the minimization of area of a three dimensional body under minimal width constraint.

\section{Support function parametrization}
\label{support}

This section recalls some of the main properties of the support function, as well as the properties which will be used in order to implement numerically the various constraints of interest in this work. The references \cite{BHcw12}, \cite{schneider},  \cite{AGcwidth} and \cite{supportfSGJ} contain more details about this subject.

Let $B$ be a convex subset of $\Bbb{R}^d$. The support function of $B$ is defined on the unit sphere $\Bbb{S}^{d-1}$ by
\[ p(\theta) = \sup_{x \in B} \  \theta \cdot  x ,\]
where the dot represents the usual Euclidean dot product. Geometrically, $p(\theta)$ represents the distance from the origin to the tangent plane $\alpha$ to $B$ such that $\alpha$ is orthogonal to $\theta$, taking into account the orientation given by $\theta$. Therefore the sum of the values of the support function for two antipodal points gives the \emph{width} or diameter of $B$ in the direction defined by these two points. This shows that bounds on the width of $B$ can be expressed by inequalities of the type
\begin{equation*}
 w \leq p(\theta)+p(-\theta) \leq W \text{ for every } \theta \in \Bbb{S}^{d-1}
 \label{diam-bound}
 \end{equation*}
and a constant width constraint can be expressed by
\begin{equation*}
 w = p(\theta)+p(-\theta) \text{ for every } \theta \in \Bbb{S}^{d-1}.
 \label{cw-constr}
\end{equation*}

As already shown in \cite{antunesMM16}, it is possible to impose inclusion constraints when dealing with support functions. Consider two convex bodies $B_1,B_2$ with support functions given by $p_1,p_2$. Then $B_1$ is included in $B_2$ if and only if $p_1(\theta)\leq p_2(\theta)$ for every $\theta \in \Bbb{S}^{d-1}$. In the case where $B_2$ is an intersection of half-spaces the inequality $p_1(\theta) \leq p_2(\theta)$ only needs to be imposed for a finite number of directions $\theta \in \Bbb{S}^{d-1}$, corresponding to the normals to the hyperplanes determining the hyperspaces. 

Each convex body in $\Bbb{R}^d$ has its own support function, however the converse is not true. The necessary assumptions for a function $p: \Bbb{S}^{d-1}\to \Bbb{R}$ to be the support function of a convex body together with precise regularity properties of the support function are discussed in \cite[Section 1.7]{schneider}. In the following, we make the assumption that the support functions are smooth enough so that first and second partial derivatives can be computed. In Section \ref{sec:theory} it is recalled that this assumption is not too restrictive, since the class of convex shapes with smooth support function is dense in the class of convex sets with respect to the Hausdorff metric. Moreover, support functions used in the numerical simulation verify this assumption.

Given a convex set $B$ and its support function $p$, a parametrization of $\partial B$ is given by
\[ \Bbb{S}^{d-1} \ni \theta \mapsto x(\theta)= p(\theta)\theta + \nabla_\tau p(\theta) \in \Bbb{R}^d,\]
where $\nabla_\tau$ represents the tangential gradient with respect to the metric in $\Bbb{S}^{d-1}$. Note that for this parametrization the normal to $\partial B$ at the point $x(\theta) \in \partial B$ is given by $\theta$. The convexity constraint can be expressed by the fact that the principal curvatures of the surface are everywhere non-negative. In the following, the aspects related to convexity are detailed and the presentation is divided with respect to the dimension. 

\subsection{Dimension 2}

In $\Bbb{R}^2$, the unit circle $\Bbb{S}^1$ is identified to the interval $[0,2\pi]$, therefore the parametrization of the boundary of the shape in terms of the support function becomes
\begin{equation} \begin{cases}
 x(\theta) = p(\theta) \cos \theta - p'(\theta) \sin \theta \\
 y(\theta) = p(\theta) \sin \theta + p'(\theta) \cos \theta.
\end{cases} 
\label{param2D}
\end{equation}
It is immediate to see that $\|(x'(\theta),y'(\theta))\| = p(\theta)+p''(\theta)$ and, as already underlined in \cite{BHcw12}, the convexity constraint in terms of the support function is equivalent to $p+p''\geq 0$, in the distributional sense.

\subsection{Dimension 3} 

In $\Bbb{R}^3$ it is classical to consider the parametrization of $\Bbb{S}^2$ given by
\begin{equation}
 \bo{n}= \bo{n}(\phi,\psi) \mapsto (\sin \phi\sin \psi,\cos \phi \sin \psi, \cos \psi),\ \phi \in [-\pi,\pi),\psi \in (0,\pi).
 \label{sphereParam}
 \end{equation}
As recalled in \cite{supportfSGJ}, if $p = p(\phi,\psi)$ is a $C^1$ support function then a parametrization of the boundary is given by
\begin{equation}
 \bo{x}_p (\phi,\psi) = p(\phi,\psi) \bo n+\frac{p_\phi(\phi,\psi)}{\sin^2 \psi} \bo n_\phi+p_\psi(\phi,\psi) \bo n_\psi
 \label{supp3Dparam}
\end{equation}
Moreover, the differential $d\bo{x}_p$ on the basis $\bo n_\phi,\bo n_\psi$ of the corresponding tangent space to $\Bbb{S}^2$ is given by
\begin{align}
  d\bo x_p \mid_{\bo n}(\bo n_\phi) & = \left(  p \sin \psi+ \frac{p_{\phi\phi}}{\sin\psi}+p_\psi \cos \psi  \right) \frac{\bo n_\phi}{\sin \psi} + \left( -\frac{p_\phi \cos \psi}{\sin \psi} + p_{\psi \phi} \right)\bo n_\psi    \notag \\ 
  d\bo x_p \mid_{\bo n}(\bo n_\psi) & = \left( \frac{p_{\phi\psi}}{\sin \psi} -\frac{p_\phi \cos \psi}{\sin^2 \psi} \right)\frac{\bo n_\phi}{\sin \psi} + (p+p_{\psi \psi} ) \bo n_\psi .
  \label{diff3D}
\end{align}
Note that $\{\bo n_\phi/\sin(\psi),\bo n_\psi\}$ is an orthonormal basis of the tangent space when $\sin\psi\neq 0$. The convexity constraint is characterized by the non-negativity of the principal curvatures. Another characterization is given by the positive definiteness of the matrix having coefficients given by the differential \eqref{diff3D} of $\bo x_p$ for every $\phi \in [-\pi,\pi)$ and $\psi \in [0,\pi)$. However, in dimension $3$ it is enough to impose a simpler condition. Indeed, if a surface has non-negative Gaussian curvature in a neighborhood of a point, then the surface is locally convex around that point. Tietze's theorem states that if a three dimensional set is locally convex around each point then it is globally convex \cite[p. 51-53]{valentineConvex}. 

As a consequence, a closed surface in dimension $3$ which has positive Gaussian curvature everywhere bounds a convex body. This is also known as Hadamard's problem \cite[p. 108]{toponogov}. Therefore, in dimension three, the convexity constraint can be imposed by assuring that the Gaussian curvature is positive at every point. More explicitly, the determinant of the matrix containing the coefficients of the differential \eqref{diff3D} is positive:   
\begin{equation}\left( p\sin \psi + \frac{p_{\phi\phi}}{\sin\psi}+p_\psi \cos \psi  \right)(p+p_{\psi \psi}) +\frac{1}{\sin \psi}\left( \frac{p_\phi \cos \psi}{\sin \psi} - p_{\psi \phi} \right)^2 > 0 
\label{convexity3D}
\end{equation} 
for every $\phi \in [-\pi,\pi),\ \forall\psi \in (0,\pi)$.

  Note that formulas \eqref{supp3Dparam}, \eqref{diff3D} and \eqref{convexity3D}  contain $\sin \psi$ in some of the denominators. In the numerical simulations the discretization of the sphere is always chosen avoiding the north and south poles of the sphere $\Bbb{S}^2$ where $\sin \psi$ cancels and where singular behavior may occur.

\section{Theoretical aspects}
\label{sec:theory}

\subsection{Existence of optimal shapes}

When dealing with shape optimization problems the existence of optimal shapes is a fundamental question. All problems studied numerically in this article deal with convex domains that are contained in a bounded set, so it is useful to define the class $\mathcal K^d$ of closed convex sets in $\Bbb{R}^d$ which are contained in a closed, large enough, ball $B$. The question of existence of solutions is greatly simplified due to the following classical result \cite[Theorem 1.8.7]{schneider}. 

\begin{thm}[Blaschke selection theorem]
	Given a sequence $\{K_n\}$ of closed convex sets contained in a bounded set, there exists a subsequence which converges to a closed convex set $K$ in the Hausdorff metric.
	\label{blaschke}
\end{thm}
For the sake of completeness, recall that the Hausdorff distance between two convex bodies $K_1,K_2$ is defined by 
\[ d_H(K_1,K_2) = \max\left\{ \sup_{x \in K_1} \inf_{y \in K_2} |x-y|,\sup_{x \in K_2} \inf_{y \in K_1} |x-y|\right\}.\]
A sequence of closed convex sets $\{K_n\}$ converges to $K$ in the Hausdorff metric if and only if $d(K,K_n) \to 0$ as $n \to \infty$.

More details regarding this result and proofs can be found in \cite[Chapter 2]{henrot-pierre-english}, \cite[Chapters 6, 7]{delfour-zolesio}. Existence results for all problems studied in this paper are a consequence of the properties listed below. These properties are classical, but are recalled below from the sake of completeness, with sketches of proof when the proofs were not readily found in the literature. If $K_1,K_2$ have support functions $p_{K_1}$ and $p_{K_2}$ then the Hausdorff distance is simply $d_H(K_1,K_2)=\|p_{K_1}-p_{K_2}\|_\infty$ \cite[Lemma 1.8.14]{schneider}.

\begin{prop}
	Convexity is preserved by the Hausdorff convergence.
\end{prop}
 For a proof see \cite[p. 35]{henrot-pierre-english}.
\begin{prop}
	If $\{K_n\}$ is a sequence of non-empty closed convex sets contained in a bounded set then the Hausdorff convergence of $K_n$ to $K$ is equivalent to the uniform convergence of the support functions $p_{K_n}$ to $p_K$ on $\Bbb{S}^{d-1}$.
	\label{prop:haussdorff-support}
\end{prop}
For a proof see \cite[Theorem 6]{salinetti-wets}. In the following, $\chi_K$ denotes the characteristic function of the set $K$. 
\begin{prop}
	Suppose that the sequence of convex sets $\{K_n\}$ converges to the convex set $K$ in the Hausdorff topology and that $K$ has non-void interior. Then $\chi_{K_n}$ converges to $\chi_K$ in $L^1$. As a consequence, $|K_n| \to |K|$ and $\mathcal H^{d-1}(\partial K_n) \to \mathcal H^{d-1} (\partial K)$ as $n \to \infty$.
	\label{prop:vol-perim}
\end{prop}

 A proof of this fact can be found in \cite[Proposition 2.4.3]{bucurbuttazzo}.


The Dirichlet-Laplace eigenvalues are solutions of 
\begin{equation}
\label{eigprob}
\left\{ \begin{array}{rcll}
-\Delta u & = & \lambda_k(\omega)u & \text{ in } \omega \\
u & = & 0 &  \text{ on } \partial \omega.
\end{array} \right.
\end{equation}
It is classical that for Lipschitz domains, the spectrum of the Dirichlet-Laplace operator consists of a sequence of eigenvalues (counted with multiplicity)
\[ 0<\lambda_1(\omega) \leq \lambda_2(\omega) \leq ... \to \infty.\]
In particular, convex sets with non-void interior enter into this framework.  Two basic properties of the Dirichlet-Laplace eigenvalues are the monotonicity with respect to inclusion and the scaling property:
\[\omega_1 \subset \omega_2 \Rightarrow \lambda_k(\omega_1) \geq \lambda_k(\omega_2) \text{   and   } \lambda_k(t \omega) = \frac{1}{t^2}\lambda_k(\omega).\]

The following property deals with the continuity of these eigenvalues with respect to the Hausdorff metric.

\begin{prop}
	If $K_n$ are convex and converge to $K$ in the Hausdorff metric then $K_n$ $\gamma$-converges to $K$ and, in particular the eigenvalues of the Dirichlet-Laplace operator are continuous: $\lambda_k(K_n) \to \lambda_k(K)$.
	\label{prop:gconv}
\end{prop} 
For more details see \cite[Theorem 2.3.17]{henroteigs}.

\begin{prop}
	Inclusion is stable for the Hausdorff convergence: $K_n \subset \Omega$, $K_n \to K$ implies $K\subset \Omega$.
	\label{prop:inclusion}
\end{prop} 
For a proof see  \cite[p. 33]{henrot-pierre-english}. 

\begin{prop}
	The diameter and width constraints are continuous with respect to the Hausdorff convergence of closed convex sets. In particular if the sequence of closed convex sets $\{K_n\}$ converges to $K$ in the Hausdorff metric and each $K_n$ is of constant width $w$ then $K$ is also of constant width $w$.
	\label{prop:width}
\end{prop}

\emph{Proof:} Property \ref{prop:haussdorff-support} recalled above shows that the Hausdorff convergence implies the uniform convergence of support functions on $\Bbb{S}^{d-1}$. Therefore, diameter and width constraints that can be expressed in pointwise form starting from the support functions are preserved, in particular, the constant width property. \hfill $\square$

In the following, $\{\phi_i\}_{i=0}^\infty$ denotes an orthogonal basis of $L^2(\Bbb{S}^{d-1})$ made of eigenvalues of the Laplace-Beltrami operator on $\Bbb{S}^{d-1}$ (the Fourier basis in 2D and the spherical harmonics in 3D). Denote by $\lambda_i\geq 0$, $i \geq 0$, the corresponding eigenvalues: $-\Delta_\tau \phi_i = \lambda_i \phi_i\text{ on } \Bbb{S}^{d-1}$. In particular $\lambda_0$ corresponds to the constant eigenfunction $\phi_0$. When studying rotors, only some particular coefficients in the spectral decomposition are non-zero. Also, in numerical approximations a truncation of the spectral decomposition is used. Therefore, it is relevant to see if such a property is preserved by the Hausdorff convergence of convex sets. Let $J\subset \Bbb{N}$ be a non-empty, possibly infinite subset of indices and denote by
{\small
\[ \mathcal F_J = \{ p : p = \sum_{i=0}^\infty \alpha_i \phi_i, p \text{ is the support function of a convex body}, \alpha_i=0  \ \ \forall i \notin J\},\]
}
i.e. convex shapes for which the coefficients of the support function in the basis $\{\phi_i\}_{i=0}^\infty$ having indices that are not in $J$ are zero. The following result holds:
\begin{prop}
	For a fixed set $J \subset \Bbb{N}$, let $\{K_n\}$ be a sequence of closed convex sets contained in a bounded set $B$ with support functions $(p_{K_n})$ contained in $\mathcal F_J$ such that $K_n$ converges to $K$ in the Hausdorff metric. Then the support function of $K$ also belongs to $\mathcal F_J$.
	\label{prop:stab_fourier}
\end{prop}

\emph{Proof:} From Property \ref{prop:haussdorff-support} it follows that the support functions of $K_n$ converge uniformly to the support function of $K$, i.e. $\|p_{K_n}-p_K\|_\infty \to 0$ as $n \to \infty$. This obviously implies the convergence in $L^2(\Bbb{S}^{d-1})$ of the support functions. In particular, if $p_{K} = \sum_{i=1}^\infty \alpha_i \phi_i$ and $p_{K_n} = \sum_{i=1}^\infty \alpha_i^n \phi_i$ then
\[ \int_{\Bbb{S}^{d-1}} (p_{K_n}-p_K)\phi_i = \alpha_i^n-\alpha_i.\]
Since the left hand side converges to zero as $n \to \infty$ for all $i$, it follows that $\alpha_i^n \to \alpha_i$ for all $i$. In particular, if $i \notin J$ then all $\alpha_i^n=0$. As a consequence $\alpha_i = 0$ for all $i \notin J$, which means that $p_K \in \mathcal F_J$. \hfill $\square$

It is possible to characterize constant width bodies by imposing that certain coefficients of the support function are zero. 

\begin{prop}
	Let $p = \sum_{i=0}^\infty \alpha_i \phi_i$ be the support function of a body with constant width. Then $\alpha_i=0$ whenever $\phi_i$ is even and non-constant on $\Bbb{S}^2$.
	\label{prop:coef-cw}
\end{prop} 

\emph{Proof:} The support function of a body of constant width satisfies $p(\theta)+p(-\theta) = w$ for every $\theta \in \Bbb{S}^2$. The conclusion follows from the fact that $\phi_i$ form an orthogonal basis of $L^2(\Bbb{S}^{d-1})$ and $\int_{\Bbb{S}^{d-1}} \phi_i = 0$ for all $i\geq 1$.\hfill $\square$


Starting from these properties, the existence of solutions for all problems studied numerically in the following section can be proved. Recall that $\lambda_k(\Omega)$ denotes the $k$-th eigenvalue of the Dirichlet-Laplace operator defined by \eqref{eigprob}.

\begin{problem} [\bf Minimizing Dirichlet-Laplace eigenvalues under convexity constraint.]
	\label{prob:eigs-convex}
	$$ \min \{\lambda_k(\Omega) : \Omega \subset \mathbb{R}^N, \Omega \text{ convex }, | \Omega | = c.\} 
	$$
\end{problem}
The existence of solutions for this problem is proved in \cite[Theorem 2.4.1]{henroteigs}. It is a direct consequence of Theorem \ref{blaschke} and Properties \ref{prop:vol-perim}, \ref{prop:gconv} above. 

\begin{problem}[\bf Minimizing Dirichlet-Laplace eigenvalues under convexity and constant-width constraints.]
	\label{prob:eigs-cw}
	\[ \min \{\lambda_k(\Omega) : \Omega \subset \Bbb{R}^N, \Omega \text{ convex with fixed constant width }, |\Omega| = c. \} \]
\end{problem}
For proving the existence in this case, choose a minimizing sequence $\{K_n\}$ which, by Theorem \ref{blaschke}, up to a subsequence, converges to $K$. In view of the Properties \ref{prop:vol-perim}, \ref{prop:gconv} and \ref{prop:width} above $K$ is indeed a solution.

\begin{problem}[\bf Minimizing the volume under constant width constraint.]
	\label{prob:minvol-cw}
	\[ \min \{ |\Omega| :  \Omega \subset \Bbb{R}^N, \Omega \text{ convex with fixed constant width }\} \]
\end{problem} 
 The existence follows from Properties \ref{prop:vol-perim}, \ref{prop:width} when working with a converging minimizing sequence, which exists by Theorem \ref{blaschke}.
 
 The next problem considered concerns rotors of minimal volume. A rotor is a convex shape that can be rotated inside a polygon (or polyhedron) while always touching every side (or face). A survey on rotors in dimension two and three can be found in \cite{goldberg_rotors}. In particular, the article \cite{goldberg_rotors} describes which coefficients are non-zero in the spectral decomposition of the support function of rotors, using Fourier series in 2D or spherical harmonics in 3D. It turns out that the earliest complete development on the subject was published in $1909$ by Meissner \cite{meissner_rotors}. More details and proofs of the claims in the papers described above can be found in \cite{groemer}.
 
 In dimension two, every regular $n$-gon admits non-circular rotors and they are characterized by the fact that only the coefficients for which the index has the form $nq\pm 1$ are non-zero, where $q$ is a positive integer. In dimension three, there are only three regular polyhedra which admit rotors: the regular tetrahedron, the cube and the regular octahedron. The rotors in a cube are bodies of constant width. For rotors in a tetrahedron the only non-zero coefficients correspond to the spherical harmonics with indices $0, 1,2$ and $5$, while in the case of the octahedron the non-zero coefficients have indices $0,1$ and $5$. The constant term in the spectral decomposition of the support function of a rotor corresponds to the inradius of the domain.

\begin{problem}[\bf Rotors of minimal volume.]
	\label{prob:rotors}
	For $P$ a polygon (polyhedron) which admits rotors, solve
	\[ \min \{ |\Omega| : \Omega\subset P, \Omega \text{ is a rotor}\}.\]
\end{problem}
Assuming $P$ admits rotors, the existence of rotors of minimal volume is guaranteed by Theorem \ref{blaschke} and properties \ref{prop:inclusion}, \ref{prop:stab_fourier}. Moreover, the fact that the Hausdorff limit of rotors is still a rotor comes from Property \ref{prop:haussdorff-support}. It is enough to choose the normal directions orthogonal to the sides of the polygon (polyhedron) and observe the limit of the corresponding support functions evaluated at these directions.
 
 The following problem was considered in \cite{oudetCW} and consists in minimizing the area under minimal width constraint.

\begin{problem}
	\label{prob:area-minw}Minimize the area of a convex three dimensional shape $\omega$ which has minimal width equal to $1$. 
\end{problem}
In \cite[Proposition 2.4.3]{bucurbuttazzo}, it is proved that of two convex bodies $A,B \subset \Bbb{R}^d$ verify $A \subset B$ then $\mathcal H^{d-1}(\partial A) \leq \mathcal H^{d-1}(\partial B)$. Every body with minimal width equal to $1$ contains three mutually orthogonal segments of length $\geq 1$. By convexity and the property above it is immediate to see that the volume of any body in $\mathcal K$ is at least $\frac{1}{8}$. Using the isoperimetric inequality, it can be seen that minimizing sequences exist and the existence of a solution comes from Theorem \ref{blaschke} and Property \ref{prop:width}.

As an application for the inclusion constraint, the Cheeger set associated to some convex domains in dimension two and three is considered. 
\begin{problem}[\bf Cheeger sets.] 
	\label{prob:Cheeger}
	The Cheeger set associated to a convex domain $\Omega \subset \Bbb{R}^d$ is the solution of the problem
	\[ \min_{X \subset \Omega} \frac{\mathcal{H}^{d-1}(\partial X)}{|X|},\]
	where the minimum is taken over all convex sets $X$ contained in $\Omega$.
\end{problem} 

The Cheeger sets are extensively studied and it is not the objective to present the subject in detail here. In dimension two there is an efficient characterization which allows the analytical computation of Cheeger sets for a large class of domains \cite{KLR06}. Computational approaches based on various methods were introduced in \cite{LROconvex}, \cite{CCP09}, \cite{cfm09} and \cite{BBF18}. Note that in dimension two, the convexity of $\Omega$ implies the convexity of the optimal Cheeger set. In dimension three this is no longer the case. However, one can prove that there exists at least one convex optimal Cheeger set \cite{LROconvex}. Existence of Cheeger sets is a consequence of Theorem \ref{blaschke} and Properties \ref{prop:vol-perim}, \ref{prop:inclusion}.

\subsection{Convergence results}

The numerical approach described in Section \ref{sec:numerics} uses a truncation of the spectral decomposition of the support function. Therefore, as underlined in \cite{BHcw12}, it is needed to prove that increasing the number of non-zero coefficients $N$ in the parametrization gives optimal shapes which converge to the solution of the original problem. 

In the following, denote by 
\[ \mathcal F_N = \{ p : p = \sum_{i=0}^N \alpha_i \phi_i, p \text{ is the support function of a convex body}\}.\]
This corresponds to the notation $\mathcal F_J$ used previously, with $J=\{0,1,2,...,N\}$. 

Denote by $\mathcal K_N^d$ the class of convex sets in $\Bbb{R}^d$ whose support functions belong to $\mathcal F_N$. In \cite[Appendix]{schneider} it is proved that for $N$ large enough, $\mathcal F_N$ is not trivial and therefore $\mathcal K_N^d$ is also non-trivial. Property \ref{prop:stab_fourier} proved in the previous section shows that $\mathcal K_N^d$ is closed in the Hausdorff metric. Therefore the existence of solutions can be shown for all problems recalled in the previous section by replacing the class of convex sets with $\mathcal F_N$.

An important question, which is not obvious at first sight, is whether a general convex body $K$ can be approximated in the Hausdorff metric by convex bodies $K_n$ with support functions in some $\mathcal F_{N_n}$. This is proved in \cite[Section 3.4]{schneider}. Namely, the following property holds:

\begin{prop}
	Let $K$ be a convex body and $\varepsilon>0$. Then there exists a positive integer $N_\varepsilon>0$ and a convex set $K_\varepsilon$ with support function in $\mathcal F_{N_\varepsilon}$ such that $d_H(K,K_\varepsilon)<\varepsilon$. 
	\label{prop:approx-finite}
\end{prop}

\begin{rem}
	Following the remarks in \cite[p. 185]{schneider}, starting from a body of constant width $K$, the smoothing procedure preserves the constant width. Moreover, the approximation of $K$ in $\mathcal K^d_{N_\varepsilon}$ is obtained by truncating the spectral decomposition of the regularized support function. This also preserves the constant width, which by Property \ref{prop:coef-cw} is equivalent to the fact that the coefficients of the even basis functions are zero. Therefore, if $K$ is of constant width in the previous proposition, its approximation $K_\varepsilon \in \mathcal F_{N_\varepsilon}$ can also be chosen of the same constant width.
	\label{approx:finite-cw}
\end{rem}

In practice, however, it is often necessary to impose some other constraints, like fixed volume, area, minimal width, etc. Below we give another variant of this property for constraints of the form $\{\mathcal C(K) \geq c\}$ where $\mathcal C$ is a continuous function with respect to the Hausdorff metric which is homogeneous of degree $\alpha>0$: $\mathcal C(\eta K) = \eta^\alpha\mathcal C(K)$ for $\eta>0$. This includes many constraints of interest, like area, perimeter, minimal width, diameter, etc.

\begin{prop}
	Let $\mathcal C:\mathcal K^n \to \Bbb{R}_+$ be a continuous functional, positively homogeneous of degree $\alpha>0$. Let $K$ be a convex body which satisfies $\mathcal C(K) \geq c$, for some fixed $c>0$ and $\varepsilon>0$. Then there exists a positive integer $N_\varepsilon>0$ and a convex set $K_\varepsilon$ with support function in $\mathcal F_{N_\varepsilon}$ such that $d_H(K,K_\varepsilon)<\varepsilon$ and $\mathcal C(K_\varepsilon) \geq c$. 
	\label{approx:finite-constr}
\end{prop}

\emph{Proof:} Property \ref{prop:approx-finite} implies the existence of a sequence $K_n$ converging to $K$ in the Hausdorff metric such that $K_n \in \mathcal K^d_{N_n}$ for some $N_n>0$. Since $\mathcal C$ is continuous it follows that $\mathcal C(K_n) \to \mathcal C(K)$. Define the new sets $K_n'$ by 
$ K_n' =
 \left(\frac{\mathcal C(K)}{\mathcal C(K_n)}\right)^{1/\alpha}K_n$.
Then obviously $\mathcal C(K_n') =\mathcal C(K)\geq c$ and 
\[ d_H(K_n',K) \leq d_H(K'_n,K_n)+d_H(K_n,K). \]
The relation between the Hausdorff distance and the support functions implies that
\[ d_H(K_n',K_n) =  \|p_{K_n'}-p_{K_n}\|_\infty = \left|1-\left( \frac{\mathcal C(K)}{\mathcal C(K_n)}\right)^{1/\alpha}\right| \|p_{K_n}\|_\infty \to 0 \text{ as } n \to \infty. \]
Therefore it is possible to approximate $K$ in the Hausdorff metric with the sequence $K_n'$ which also verifies the constraint $\mathcal C(K_n') \geq c$.\hfill $\square$

Now we are ready to prove the following approximation result.

\begin{thm}
	Let $\mathcal G$ be a continuous functional defined on the class of closed convex sets. Consider a constraint function $\mathcal C$ which is continuous for the Hausdorff metric and let $c>0$. In the following $\mathcal M$ denotes one of the following: $\mathcal K_N^d$, $\{K \in \mathcal K_N^d : \mathcal C(K) \geq c\}$, or the set of sets in $\mathcal K_N^d$ of fixed constant width $w$. Denote by $K_N$ a solution of 
	\begin{equation} \min_{K \in \mathcal K_N^d\cap \mathcal M} \mathcal G(K).
	\label{pb:finite}
	\end{equation}
	Then any converging subsequence of $K_N$ converges in the Hausdorff metric to a solution $K$ of 
	\begin{equation} \min_{K \in \mathcal K^d\cap \mathcal M} \mathcal G(K).
	\label{pb:continuous}
	\end{equation}
	\label{thm:approx}
\end{thm}

\emph{Proof:} First, let us note that the existence of solutions to problems \eqref{pb:finite} and \eqref{pb:continuous} is immediate using Theorem \ref{blaschke}, Property \ref{prop:stab_fourier} the fact that $\mathcal M$ is closed and the continuity of $\mathcal G$.

Denote by $K$ a solution of \eqref{pb:continuous}. By the results shown in Properties \ref{prop:approx-finite},  \ref{approx:finite-constr} and Remark \ref{approx:finite-cw} there exists a sequence of convex sets $L_{N_n} \in \mathcal K^d_{N_n}\cap \mathcal M$ such that $L_{N_n} \to K$ in the Hausdorff metric. 

In the following denote by $K_N$ a solution of \eqref{pb:finite} for $N \geq 1$. It is obvious that $\mathcal K^d_{N_1} \subset \mathcal K^d_{N_2}\subset \mathcal K^d$ for $N_1 \leq N_2$. As an immediate consequence  $\mathcal G(K_{N_1}) \geq \mathcal G(K_{N_2}) \geq \mathcal G(K)$ for $N_1 \leq N_2$. Therefore, the sequence $\{\mathcal G(K_N)\}_{N \geq 1}$ is non-increasing and bounded from below, which implies the existence of the limit $\lim_{N\to \infty} \mathcal G(K_N) = \ell\geq \mathcal G(K)$. By the optimality of $K_N$ we have $\mathcal G(L_{N_n}) \geq \mathcal G(K_{N_n})$ and since $d_H(K,L_{K_n})\to 0$, by the continuity of $\mathcal G$ we obtain
\[ \mathcal G(K) = \lim_{n \to \infty} \mathcal G(L_{N_n}) \geq \lim_{n\to \infty} \mathcal G(K_{N_n}) = \ell \geq \mathcal G(K). \]
The inequality above implies that $\ell = \mathcal G(K)$. Furthermore, if a subsequence of $\{K_N\}_{N \geq 1}$ converges to $K'$ in the Hausdorff metric it follows that $K' \in \mathcal M$ and $\mathcal G(K') = \ell = \min_{K \in \mathcal K^d\cap \mathcal M} \mathcal G(K)$. Therefore, every limit point for $\{K_N\}_{N \geq 1}$ in the Hausdorff metric is a minimizer of \eqref{pb:continuous}. 

\hfill $\square$

Theorem \ref{thm:approx} motivates our numerical approach. In order to obtain an approximation of solutions of the shape optimization problems considered, a truncation of the spectral decomposition is used. The theoretical result states that the solutions of the finite dimensional minimization problems obtained converge to the solution of the original problem.

\section{Numerical framework}
\label{sec:numerics}

In numerical shape optimization, shapes are represented using a finite number of parameters. Theorem \ref{thm:approx} shows that considering as parameters the coefficients of the truncation of a spectral decomposition is appropriate, since solutions of the resulting finite dimensional optimization problems converge to solutions of the convex shape optimization problem. This type of numerical approaches was already used in other contexts in \cite{BHcw12}, \cite{antunesMM16}, \cite{antunesf-vol}, \cite{braxtonformulas}, \cite{antunesf-per}. Using such systems of orthogonal basis representations has further advantages which will be underlined below. Again, for the clarity of exposition, we divide the presentation following the dimension.

\subsection{Dimension 2}
\label{sec:dim2}
 We approximate the support function by a truncated Fourier series
\begin{equation}
 p(\theta) = a_0 + \sum_{k=1}^N \left( a_k \cos k\theta + b_k \sin k\theta  
\right)
\label{fouDecomp}
\end{equation}
As stated in Section \ref{support}, in order for $p$ to be the support function of a convex set in $\Bbb{R}^2$ we need to have $p''(\theta)+p(\theta) \geq 0$ for every $\theta \in [0,2\pi)$. In \cite{BHcw12} the authors provide an exact characterization of this condition in terms of the Fourier coefficients, involving concepts from semidefinite programming. In \cite{antunesMM16} the author provides a discrete alternative of the convexity inequality which has the advantage of being linear in terms of the Fourier coefficients. We choose $\theta_m = 2\pi m/M_c,\ m=1,2,...,M_c$ for some positive integer $M_c$ and we impose the inequalities $p(\theta_m)+p''(\theta_m) \geq 0$ for $m=1,...,M_c$. As already shown in \cite{antunesMM16} we obtain the following system of linear inequalities
\begin{equation}
 \begin{pmatrix}
 1 & \alpha_{1,2} & \cdots & \alpha_{1,N} & \beta_{1,2} & \cdots & \beta_{1,N} \\
 \vdots & \vdots  & \ddots & \vdots & \vdots & \ddots & \vdots \\
 1 & \alpha_{M_c,2} & \cdots & \alpha_{M_c,N} & \beta_{M_c,2} & \cdots & \beta_{M_c,N}
 \end{pmatrix}
 \begin{pmatrix}
 a_0 \\ a_2 \\ \vdots \\ a_N \\ b_2 \\ \vdots \\ b_N
 \end{pmatrix}
 \geq
 \begin{pmatrix}
 0\\ \vdots \\ 0
 \end{pmatrix}
\label{eq:conv-lin-ineq}
\end{equation}
where $\alpha_{m,n}=(1-n^2)\cos(n\theta_m)$ and $\beta_{m,n} = (1-n^2)\sin(n\theta_m)$.

Next we turn to the constant width condition $p(\theta)+p(\theta+\pi)=w$ for every $\theta \in [0,2\pi)$, which is equivalent to $a_0=w/2$ and $a_{2k}=b_{2k} = 0,\ k=1,...,N$. This was already noted in \cite{BHcw12}. 

An upper bound $W$ on diameter can be introduced as a constraint for the support function 
as follows
\[  p(\theta)+p(\theta+\pi) \leq W, \theta \in [0,2\pi).\]
In the computations we consider a discrete version of the above inequality. Pick $\theta_m = 2\pi m/M_d$, $m=1,2,...,M_d$ for some positive integer $M_d$ and impose the following linear inequalities
\begin{equation*}
 p(\theta_m)+p(\theta_m+\pi) \leq W,\ m=1,...,M_d.
 \label{diam-bounds-discrete}
 \end{equation*}
In order to impose a lower bound on the diameter it is enough to pick one direction $\theta$ and use the constraint 
\[ p(\theta) + p(\theta+\pi) \geq w.\]
It is also possible to consider variable lower and upper bounds on the width of the body which depend on $\theta$.

Let us now recall the formulas for the area and perimeter of a two dimensional shape in terms of the Fourier coefficients of the support function. The perimeter is simply equal to $P(p) = 2\pi a_0$, which is linear in terms of the Fourier coefficients. As already stated in \cite{BHcw12} the area of a convex shape having support function $p$ with the Fourier decomposition \eqref{fouDecomp} is given by
\begin{equation}
 A(p) = \pi a_0^2 +\frac{\pi}{2} \sum_{i=1}^N  (1-k^2)(a_k^2 + b_k^2).
 \label{eq:area-2d}
 \end{equation}
Note that $a_1$ and $b_1$ do not contribute to the area computations as modifying $a_1,b_1$ only leads to translations of the shape defined by $p$.

\subsection{Dimension 3}
\label{sec:dim3}
In \cite{antunesf-per} the authors parametrized three dimensional domains by their radial function using spherical harmonics. In our case we parametrize the support function using a finite number of spherical harmonics
\begin{equation}
 p(\phi,\psi) = \sum_{l=0}^N \sum_{m=-l}^l a_{l,m} Y_l^m(\psi,\phi) \label{sphDecomp}
 \end{equation}
for a given positive integer $N$. The spherical harmonics are defined by
\[ Y_l^m(\psi,\phi) = \begin{cases}
  \sqrt{2} C_l^m \cos(m\phi) P_l^m (\cos \psi) & \text{ if }m >0\\
  C_l^0 P_l^0(\cos \psi) & \text{ if }m=0\\
  \sqrt{2} C_l^m \sin(-m\phi) P_l^{-m} (\cos \psi) & \text{ if }m<0, 
\end{cases}\]
where $P_l^m$ are the associated Legendre polynomials and 
\[ C_l^m = \sqrt{ \frac{(2l+1)(l-|m|)!}{4\pi(l+|m|)!}} \]
are normalization constants. 

The convexity constraint is imposed by considering a discrete version of \eqref{convexity3D}. We construct a family of $M_c$ evenly distributed points on the unit sphere, for example like described in \cite[Section 3]{antunes3D}. We denote by $(\phi_i,\psi_i)$ $i=1,...,M_c$ the corresponding pairs of angles given by the parametrization \eqref{sphereParam}. We impose that the convexity condition \eqref{convexity3D} is satisfied at points given by $(\phi_i,\psi_i),\ i=1,...,M_c$. As in the two dimensional case, width inequality constraints can be handled in a similar way, by imposing inequalities of the type
\[ w_i \leq p(\theta_i)+p(-\theta_i) \leq W_i\]
at points $\theta_i=\bo (\phi_i,\psi_i)$ (see \eqref{sphereParam}).

The constant width condition is $p(\theta)+p(-\theta) = w$ for every $\theta \in \Bbb{S}^2$. As shown in Property \ref{prop:coef-cw} this amounts to considering only odd spherical harmonics in the decomposition \eqref{sphDecomp}, except for the constant term. This corresponds to spherical harmonics $Y_l^m$ for which the index $l$ is odd. In the following, denote by $h$ the part of the support function containing the non-constant terms:  $h = p-\frac{1}{4\pi}\int_{\Bbb{S}^2}pd\sigma$.

The area and volume of a convex body of constant width $w$ in dimension three can be computed explicitly in terms of the spherical harmonics coefficients. Indeed, in \cite[Theorem 2]{AGcwidth}, the following formulas are provided:
\begin{equation}
 V = \frac{\pi}{6} w^3 - \frac{w}{2} \mathcal{E}(h)
 \label{eq:vol3D}
 \end{equation}
\begin{equation}
 A =  \pi w^2 - \mathcal{E}(h).
\label{eq:area3D}
\end{equation}
where $\displaystyle \mathcal{E}(h) = \int_{\Bbb{S}^{2}} \left( \frac{1}{2} |\nabla_\tau h|^2 - h^2\right) d A$. 
The formulas in \cite{AGcwidth} are for a body of constant width $2w$, which we transform so that they correspond to a body of width $w$. The spherical harmonics $Y_l^m$ form an orthonormal family, therefore it follows that $\mathcal{E}(p)$ can be computed explicitly in terms of the coefficients $a_{l,m}$ and the eigenvalues $\lambda_{l,m}$ corresponding to the spherical harmonics $Y_{l,m}$:
\begin{equation}
 \mathcal{E}(h) = \sum_{l=1}^N \sum_{m=-l}^l \left(\frac{\lambda_{l,m}}{2}-1\right)a_{l,m}^2. 
 \label{eq:energy}
\end{equation}
As a consequence, when dealing with bodies of constant width, the volume and the area have explicit formulas in terms of the coefficients $a_{l,m}$ of the decomposition \eqref{sphDecomp}.  

We note that it is also possible to compute explicitly the area of a general convex body, using the coefficients of the support function. Indeed, Lemma 1 from \cite[Section 5]{AGcwidth} is valid for general support functions $h$, not only those corresponding to a constant width body. Therefore the area of a convex body $B$ is also given by \eqref{eq:area3D}, where $w = 2a_{0,0}Y_0^0$. Also following the results stated in \cite{AGcwidth} it should also be possible to compute the volume explicitly using the Gaunt coefficients involving integrals on the sphere of products of three spherical harmonics. In our computations, for general bodies parametrized using their support function, we use the divergence theorem. The volume of a convex body $\omega$ is computed as the integral on $\partial \omega$ of a vector field $V$ with divergence equal to one. For simplicity choose $V = \frac{1}{3}\bo x = \frac{1}{3}(x,y,z)$ and we integrate $V\cdot n$ on $\partial \omega$. It is straightforward from the parametrization \eqref{supp3Dparam} that $\bo x(\theta)\cdot n(\theta) = p(\theta)$ for $\theta \in \Bbb{S}^2$.

\subsection{Visualization of results}
 We briefly present how the results are visualized. The variables in the optimization algorithm, and therefore, the output obtained are coefficients of a truncated spectral decomposition of the support function. Given such a family of Fourier coefficients (spherical harmonics coefficients) it is possible to evaluate the support function and its derivatives at any point in the unit circle (unit sphere in dimension three). 
 
 Once the values of the support function and its derivatives are known at a family of discretization points it is possible to use the formula \eqref{param2D} in dimension two (formula \eqref{supp3Dparam}  in dimension three) in order to find the associated points on the boundary of the domain. When such a family of points is known a simple contour plot is made in dimension two. 
 In dimension three, the Matlab command \texttt{convhull} is used to generate a triangulation of the surface bounding the convex body. Then the command \texttt{patch} is used to plot the surface of the convex body.

\section{Computation of the Dirichlet-Laplace eigenvalues}
\label{sec:fsol}

The Dirichlet-Laplace eigenvalue problem \eqref{eigprob} is solved numerically using the Method of Fundamental Solutions (MFS)~\cite{MR1849238,MR3070537}. This method does not need the construction of a mesh, is precise and has low computational time in dimensions two and three. We consider a fundamental solution of the Helmholtz equation,
\begin{equation}
\Phi_{\lambda}(x)=\frac{i}{4}H_{0}^{(1)}(\sqrt{\lambda}\left|
x\right|)\end{equation} 
and
\begin{equation} \Phi_{\lambda}(x)=\frac{e^{i\sqrt{\lambda}\left|
x\right| }}{4\pi\left| x\right|},\end{equation}
respectively for 2D and 3D cases, where $H_{0}^{(1)}$ denotes the first Hankel function. For a fixed value of $\lambda,$ the MFS approximation is a linear combination
\begin{equation}
\label{chaomfs}
\sum_{j=1}^{m}\alpha_{j}\Phi_{\lambda}(\cdot-y_{j}),\end{equation}
where the source points $y_j$ are placed on an admissible source set, for instance the boundary of a bounded open set
$\hat{\omega}\supset\bar{\omega}$, with $\partial\hat{\omega}$ surrounding
$\partial\omega$. By construction, the MFS approximation satisfies the PDE of the eigenvalue problem \eqref{eigprob} and we can just focus on the approximation of the boundary conditions, which can be justified by density results (see, for example \cite{MR3070537}).

Next, we give a brief description of the numerical procedure for calculating the Dirichlet-Laplace eigenvalues. We define two sets of points $W=\left\{w_i,\ i=1,...,n\right\}$ 
and $X=\left\{x_i,\ i=1,...,m\right\}$, almost uniformly distributed on the boundary $\partial\omega$, with $n<m$ and the set of source points,
$Y=\left\{w_i+\alpha n_i,,\ i=1,...,n\right\}$ 
where $\alpha$ is a positive parameter and $n_i$ is the unitary outward normal vector at the point $w_i$. We consider also some interior points $z_i$,  $i=1,...,p$ with $(p<m)$ randomly chosen in $\omega$ and used the Betcke-Trefethen subspace angle~\cite{MR2178637}. After defining the matrices
\begin{equation}
\mathbf{M}_1(\lambda)=\left[\Phi_{\lambda}(x_{i}-y_{j})\right]_{m\times
n},\end{equation} 
\begin{equation}
\mathbf{M}_2(\lambda)=\left[\Phi_{\lambda}(z_{i}-y_{j})\right]_{p\times
n}\end{equation} 
and $\mathbf{A}(\lambda)=\left[ \begin{array}{c}
\mathbf{M}_1(\lambda)  \\
\mathbf{M}_2(\lambda)
\end{array}\right]$
we compute the QR factorization
\[\mathbf{A}(\lambda)=\left[ \begin{array}{c}
\mathbf{Q}_1(\lambda)  \\
\mathbf{Q}_2(\lambda)
\end{array}\right]\mathbf{R}\]
and calculate the smallest singular value of  the block $\mathbf{Q}_1(\lambda)$, which will be denoted by $\sigma_1(\lambda)$. The approximations for the Dirichlet-Laplace eigenvalues are the local minima $\lambda$, for which $\sigma_1(\lambda)\approx0$.

\subsection{Shape derivatives for the Dirichlet Laplace eigenvalues}
\label{sec:shderiv}

Functionals like volume or area have explicit formulas in terms of the coefficients in the Fourier or spherical harmonics decomposition. This gives straightforward formulas for gradients and Hessians which can be used in optimization algorithms. When the shape functional is more complex, direct formulas are not available. Below, we present how the Hadamard shape derivatives can be used to obtain partial derivatives in terms of coefficients of the parametrization.

The Hadamard shape derivative formula shows how a shape functional $\mathcal F(\omega)$ varies when considering some perturbation of the boundary given by a vector field $V$. One way to define this to consider the derivative of the functional $t \mapsto \mathcal F((\text{Id} + tV)(\omega))$ at $t=0$. Under mild regularity assumptions it can be proved that the shape derivative may be written as a linear functional depending on the normal component of $V$. For more details one could consult \cite[Chapter 5]{henrot-pierre-english} or \cite[Chapter 9]{delfour-zolesio}. In particular, in \cite[Theorem 5.7.4]{henrot-pierre-english}, for the case of the eigenvalues of the Dirichlet-Laplace operator \eqref{eigprob} the shape derivative is
\[ \lambda_k'(\omega)(V) = -\int_{\partial \omega} \left(\frac{\partial u_k}{\partial n}\right)^2 V\cdot n d\sigma,\] 
as soon as eigenvalue $\lambda_k(\omega)$ is simple and the eigenfunction $u_k$ is in $H^2(\omega)$. This is true in the particular case of convex sets. In the following, we suppose that the functional $\mathcal{F}(\omega)$ has a Hadamard shape derivative which can be written in the form 
\begin{equation}
\mathcal{F}'(\omega)(dV) = \int_{\partial \omega} f V\cdot n d\sigma.
\label{HadamardGeneral}
\end{equation}

\subsubsection{Dimension 2} 
As already recalled in Section \ref{sec:numerics} a parametrization of the boundary of the convex shape defined by the support function $p$ is given by
\[ \begin{cases}
x(\theta) = p(\theta) \cos \theta - p'(\theta) \sin \theta \\
y(\theta) = p(\theta) \sin \theta + p'(\theta) \cos \theta.
\end{cases}\]
and a straightforward computation gives
\[ 
\begin{cases}
x'(\theta) = -(p''(\theta)+p(\theta))\sin \theta \\
y'(\theta) = (p''(\theta)+p(\theta))\cos \theta.
\end{cases}
\]
Therefore the norm of the velocity vector is given by $\|(x'(\theta),y'(\theta))\| = p''(\theta)+p(\theta)$, which helps compute the Jacobian when changing variables. Moreover, the normal to the point corresponding to parameter $\theta$ is simply $(\cos \theta,\sin \theta)$. In the following, $(a_k)_{k\geq 0}$ and $(b_k)_{k\geq 1}$ denote the Fourier coefficients of the support function: $p(\theta) = a_0 +\sum_{k\geq 1} (a_k\cos(k\theta)+b_k\sin(k\theta))$.


In order to compute the partial derivatives of the functional with respect to the Fourier coefficients it is enough to transform the perturbation of the support function into a perturbation of the boundary and use the Hadamard formula. We summarize the derivative formulas below. For simplicity, we use the abuse of notation $f(\theta) = f(x(\theta),y(\theta))$.
\begin{enumerate}
	\item \emph{Derivative with respect to $a_0$}. The corresponding boundary perturbation is $V = (\cos \theta, \sin \theta)$ and the normal component is $V\cdot n=1$. Therefore the derivative is
	\[ \frac{\partial \mathcal{F}}{\partial a_0} = \int_{\partial \omega} fd\sigma = \int_0^{2\pi} f(\theta) (p(\theta)''+p(\theta)) d\theta. \]
	\item \emph{Derivative with respect to $a_k$}. The corresponding boundary perturbation is \[V = (\cos(k\theta)\cos \theta + k\sin(k\theta)\sin \theta, \cos(k\theta)\sin \theta-k\sin (k\theta) \cos \theta)\] and the normal component is $V\cdot n=\cos (k\theta)$. Therefore the derivative is
	\[ \frac{\partial \mathcal{F}}{\partial a_k} = \int_{\partial \omega} f \cos(k\theta) d\sigma = \int_0^{2\pi} f(\theta)\cos(k\theta) (p(\theta)''+p(\theta)) d\theta. \]
	\item \emph{Derivative with respect to $b_k$}. The corresponding boundary perturbation is \[V = (\sin(k \theta)\cos \theta-k\cos(k\theta) \sin \theta, \sin(k\theta)\sin \theta+k\cos(k\theta)\cos\theta)\] and the normal component is $V\cdot n=\sin(k\theta)$. Therefore the derivative is
	\[ \frac{\partial \mathcal{F}}{\partial b_k} = \int_{\partial \omega} f \sin(k\theta)d\sigma = \int_0^{2\pi} f(\theta) \sin(k\theta)(p(\theta)''+p(\theta)) d\theta. \]
\end{enumerate}

\subsubsection{Dimension 3}

We differentiate now a functional $\mathcal{F}(\omega)$ for 3D shapes parametrized using the coefficients of the spherical harmonic decomposition \eqref{sphDecomp} of the support function. Given a general perturbation of the support function $p \mapsto p+Y$, in view of \eqref{supp3Dparam}, we find that the boundary perturbation has the form
\[V = Y \bo n + \mathcal P(\theta,\phi) \bo n_\phi + \mathcal{Q}(\theta,\phi) \bo n_\psi.\]
Since the vectors $\bo n, \bo n_\phi$ and $\bo n_\psi$ are orthogonal, the normal component is given by $V\cdot n = V\cdot \bo n = Y$. Then, using the general Hadamard derivative formula \eqref{HadamardGeneral} we find the partial derivatives of the objective function with respect to the coefficients in  \eqref{sphDecomp}:
\[ \frac{\partial \mathcal{F}}{\partial a_{l,m}} = \int_{\partial \omega} f Y_{l,m} d\sigma = \int_{-\pi}^\pi \int_0^\pi f(\phi,\psi) Y_{l,m}(\phi,\psi)\text{Jac}(\psi,\phi) d\psi d \phi,  \]
where $\text{Jac}(\psi,\phi)$ is the Jacobian function given by \eqref{convexity3D}. Indeed, the Jacobian this surface integral is computed by $\text{Jac}(\phi,\psi) = \| \partial_\phi\bo x_p \times \partial_\psi\bo x_p \|$. Note that the vectors $\partial_\phi \bo x_p, \partial_\psi \bo x_p$ are orthogonal to the normal $\bo n$ to the surface. Therefore the Jacobian reduces to $\text{Jac}(\phi,\psi) = \bo n \cdot  (\partial_\phi\bo x\times \partial_\psi \bo x)$ and using the expressions of the differential of $\bo x$ in the tangent plane given by \eqref{diff3D} we can conclude that $\text{Jac}(\phi,\psi)$ is indeed given by \eqref{convexity3D}.

\section{Applications}
\label{sec:applic}

This section shows how the proposed numerical framework applies to the various problems presented in Section \ref{sec:theory}. For each problem the shapes are discretized using a truncated spectral decomposition of the support function, as shown in Section \ref{sec:numerics}. The corresponding shape optimization problem is then approximated by a finite dimensional constrained optimization problem.

The Matlab \texttt{fmincon} function with the \texttt{interior-point} algorithm is used in each of the various problems shown below. As shown in previous sections, all constraints on the shapes are transformed into algebraic constraints on the coefficients of a spectral decomposition. The full functionality of \texttt{fmincon} is used in order to handle: linear equality or inequality constraints and non-linear constraints. The gradient of the functional and the gradients of the constraints are computed and are used in the algorithm. When possible, the Hessian matrix is also computed, and in all the other computations a LBFGS approximation is used. 

{The optimization toolbox described previously is efficient when the gradients of the objective function and of the constraints are provided. The computation time mostly depends on the size of the problem (number of variables and constraints) and on the cost of the evaluation of the objective function and its gradient. The cost of one objective function evaluation is as follows:
\begin{itemize}
	\item geometric functions (area, volume, surface area): explicit formulas, fast evaluation of the objective function.
	\item Dirichlet Laplace eigenvalues using MFS: roughly 30 seconds in dimension two and one minute in dimension three. The cost of computing higher eigenvalues is larger than for lower ones. 
\end{itemize} 
The cost of the optimization algorithm on a personal computer (Intel i7 processor, 4.2Ghz, 32GB RAM) ranges from under 10 minutes for functionals involving only geometric quantities to a couple of hours for functionals involving the Dirichlet Laplace eigenvalues in dimension three.}

\subsection{Minimize the Dirichlet-Laplace eigenvalues under volume and convexity constraints}
\label{sec:eigs}

This section presents numerical approximations of solutions to Problem \ref{prob:eigs-convex}. The theoretical and numerical study of minimization problems of the form 
\[ \min_{\omega \in \mathcal A} \mathcal \lambda_k(\omega)  \]
gained a lot of interest in the recent years. Various problems were considered, like the optimization of eigenvalues under volume constraint \cite{bucur-mink}, \cite{maz_prat}, the optimization under perimeter constraint \cite{deveper} and recently, the minimization under diameter constraint \cite{BHL17}. For many of the problems considered, explicit solutions are not known, therefore various works, like \cite{oudeteigs}, \cite{antunesf-vol}, \cite{antunesf-per} deal with the optimization of the eigenvalues for volume and perimeter constraints. Such constraints can naturally be incorporated in the functional, in view of the behaviour of the eigenvalue with respect to scaling, and therefore unconstrained optimization algorithms based on information given by the shape derivative are successfully used in practice. 

The challenges encountered when adding the convexity constraint are underlined in the study of the second eigenvalue:
\begin{equation}
 \min_{ |\omega|=1 ,\ \omega \text{ convex } } \lambda_2(\omega).
 \label{lam2convex}
 \end{equation} 
This problem is studied in \cite{henrot_oudet} where it is shown that the optimizer is not the convex-hull of two tangent disks, as conjectured before. Moreover, the boundary of the optimal set cannot contain arcs of circles. An algorithm for finding numerically the minimizer of \eqref{lam2convex} was proposed, using a penalization of the difference between the volume of the shape and the volume of its convex hull. A more precise, parametric search for the minimum of \eqref{lam2convex} was done in \cite{antunes_henrot}, giving an optimal numerical value of $\lambda_2(\omega) = 37.987$. 
 
\begin{figure}
 \centering
  \begin{tabular}{ccc}
 \includegraphics[height=0.22\textwidth]{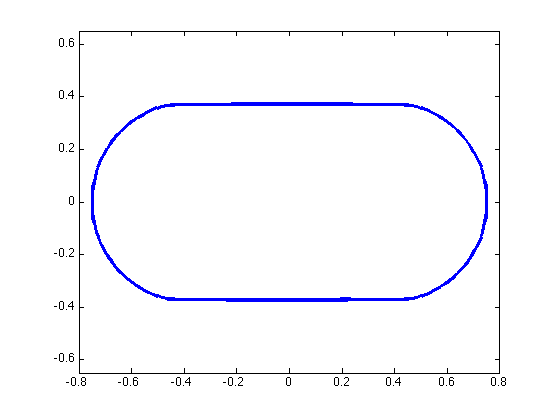}&
 \includegraphics[height=0.22\textwidth]{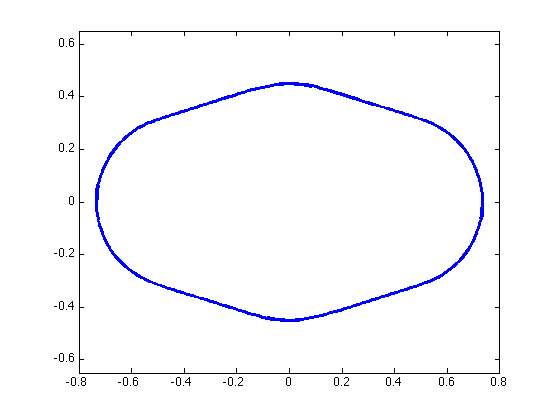}&
  \includegraphics[height=0.22\textwidth]{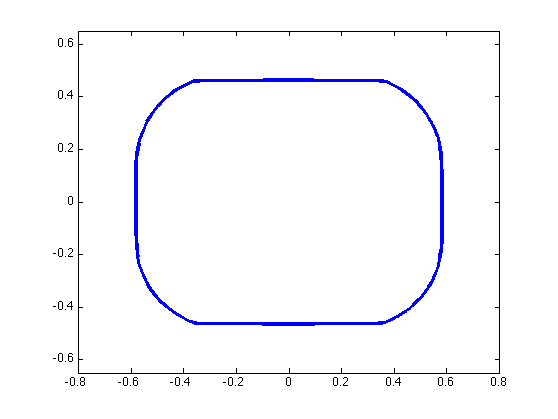}\\
    $\lambda_2=38.00$ & $\lambda_4=65.28$ & $\lambda_5=79.70$ \\
   \includegraphics[height=0.22\textwidth]{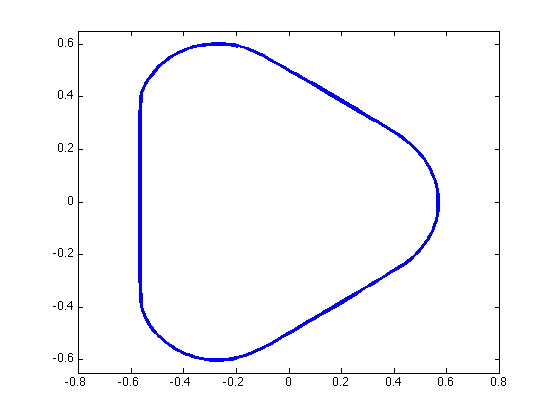}&
   \includegraphics[height=0.22\textwidth]{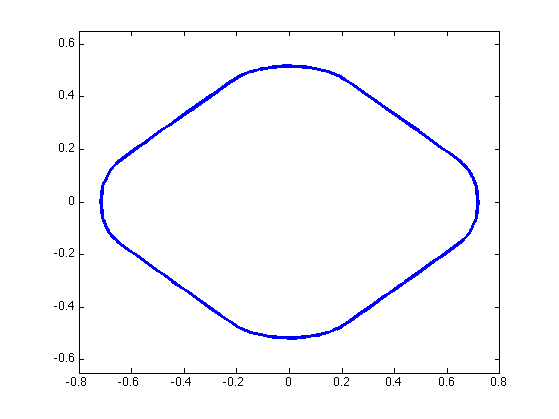}&
    \includegraphics[height=0.22\textwidth]{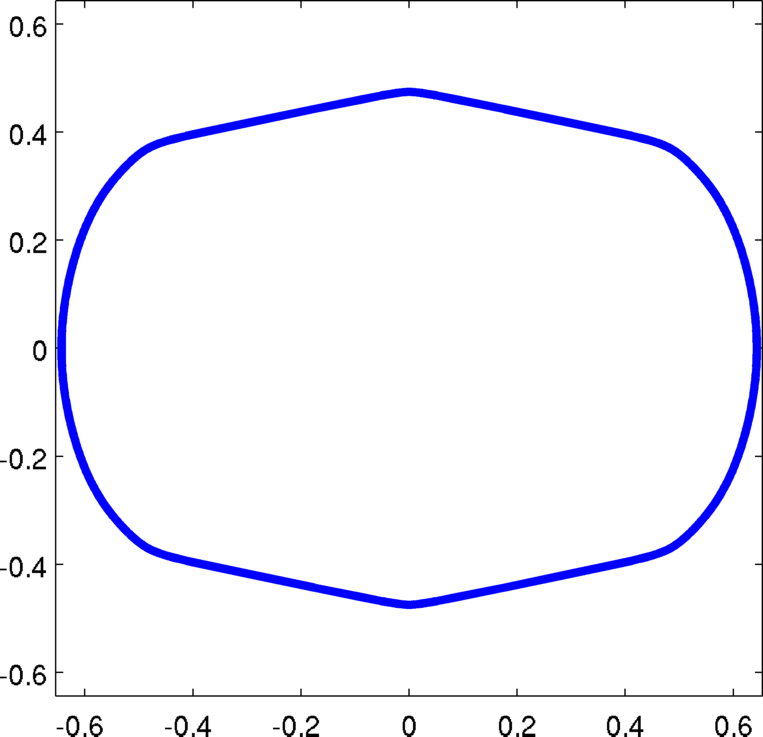}\\
        $\lambda_6=88.54 $ & $\lambda_7=109.44$ & $\lambda_8=120.58$ \\
     \includegraphics[height=0.22\textwidth]{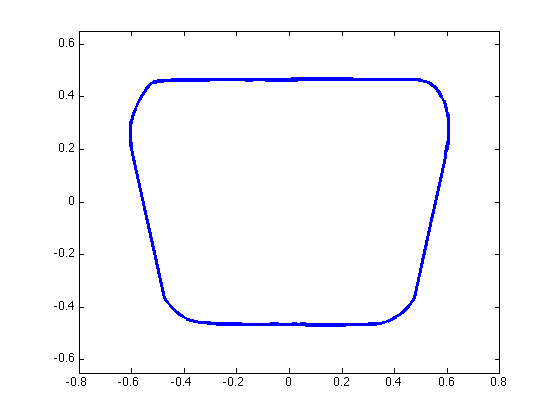}&
     \includegraphics[height=0.22\textwidth]{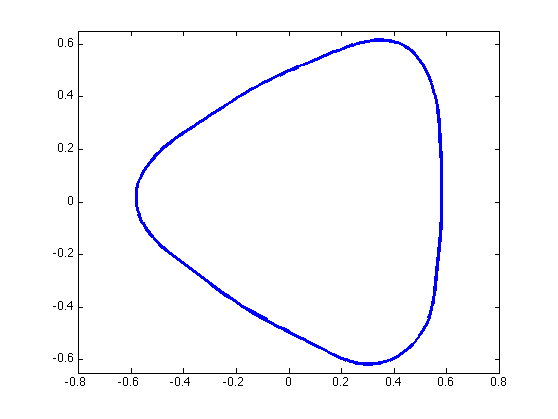}&\\
             $\lambda_9=137.38 $ & $\lambda_{10}=143.15$ &  \\
 \end{tabular}
 \caption{Minimization of the eigenvalues of the Dirichlet-Laplace operator under convexity and volume constraints in dimension two.  The numerical minimizer of the third eigenvalue is the disk, even without imposing the convexity constraint.}
 \label{fig:convexEig2D}
 \end{figure}

 The numerical algorithm proposed in this article allows us to study problem \eqref{lam2convex} in dimensions two and three. The computation of the eigenvalues is done using the method of fundamental solutions described in Section \ref{sec:fsol}. The partial derivatives with respect to the Fourier coefficients in the parametrization are computed using results shown in Section \ref{sec:shderiv} and the convexity constraint is imposed using \eqref{eq:conv-lin-ineq}.  Two dimensional results are summarized in Figure \ref{fig:convexEig2D} and it can be noted that segments are sometimes present in the boundaries of the numerical minimizers. For $k=3$ we find that the minimizer is the disk, which is in accord with the simulations performed using only a volume constraint in \cite{oudeteigs},\cite{antunesf-vol}.
 The values presented in Figure \ref{fig:convexEig2D} are obtained by rounding up the numerical optimal values and are thus upper bounds for the optimal values. The numerical simulations are made for $N=300$ ($601$ Fourier coefficients) and $M_c=5000$ points where convexity constraints are imposed.
 
\begin{figure}
 \centering
 \begin{tabular}{cccc}
 \includegraphics[width=0.2\textwidth]{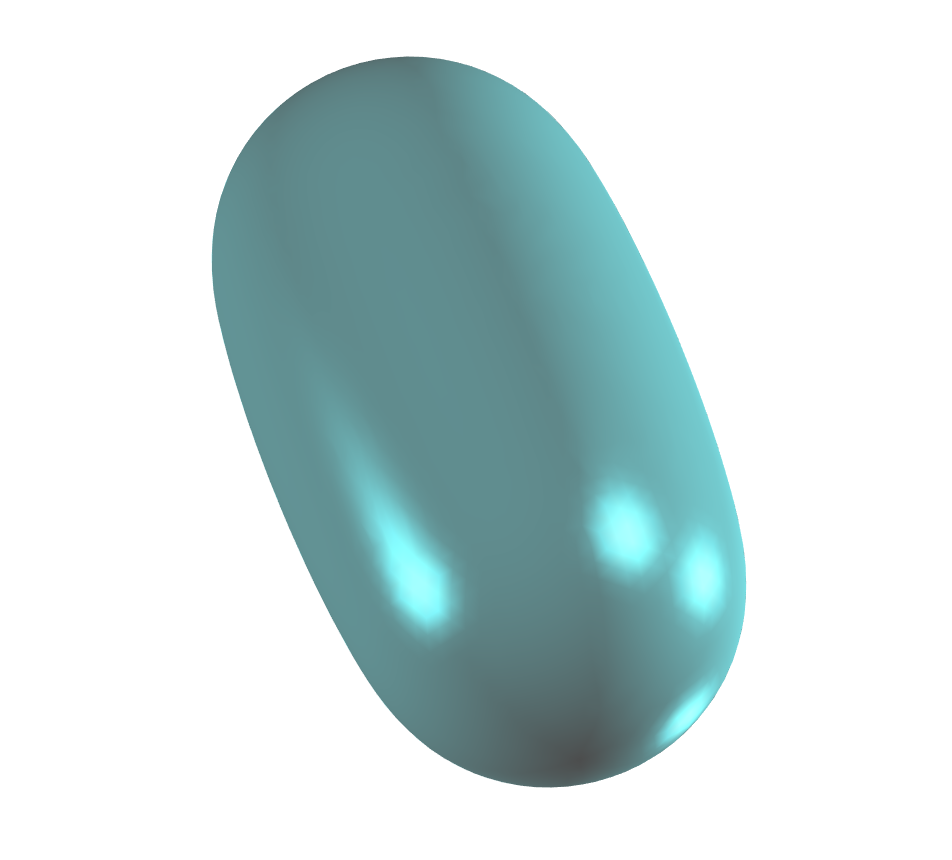}&
 \includegraphics[width=0.2\textwidth]{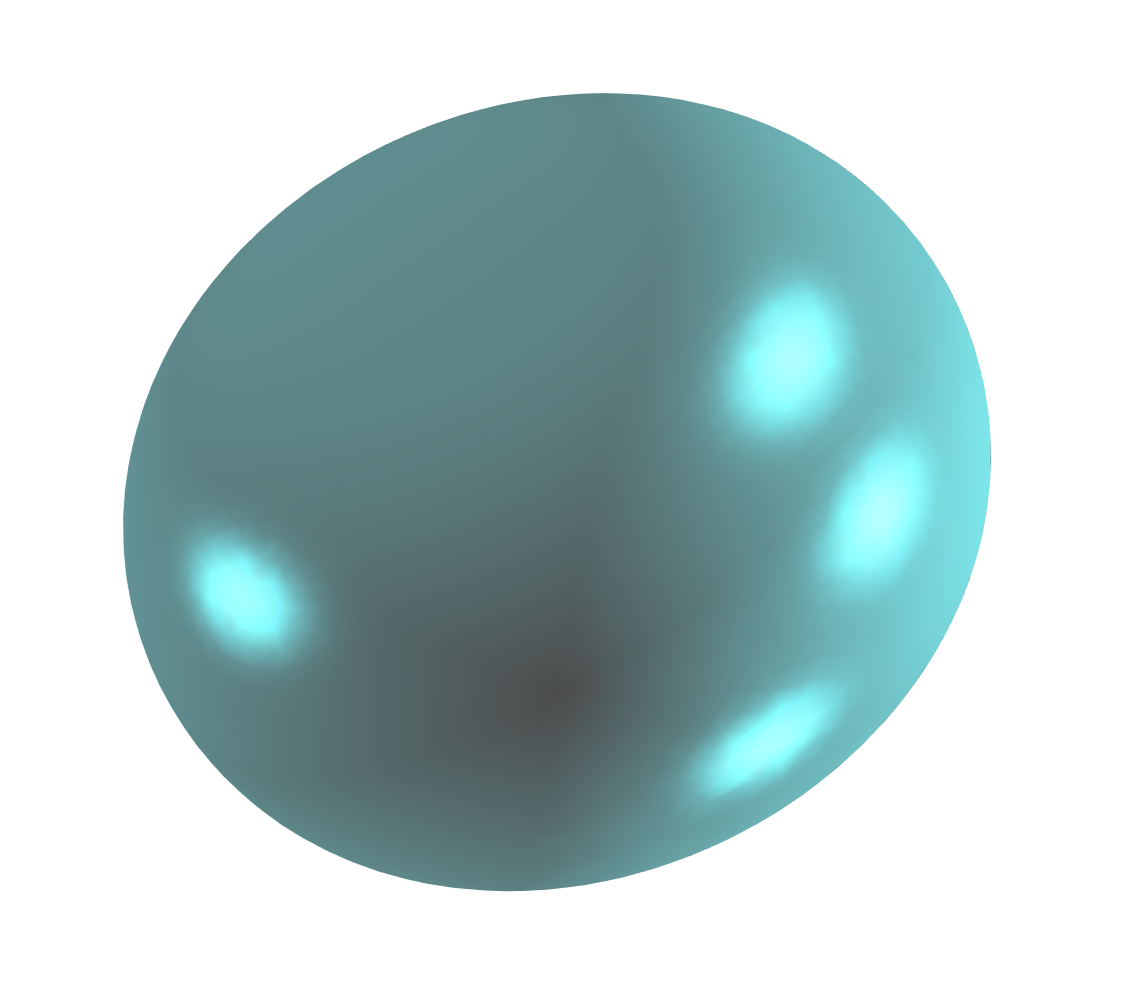}& 
 \includegraphics[width=0.2\textwidth]{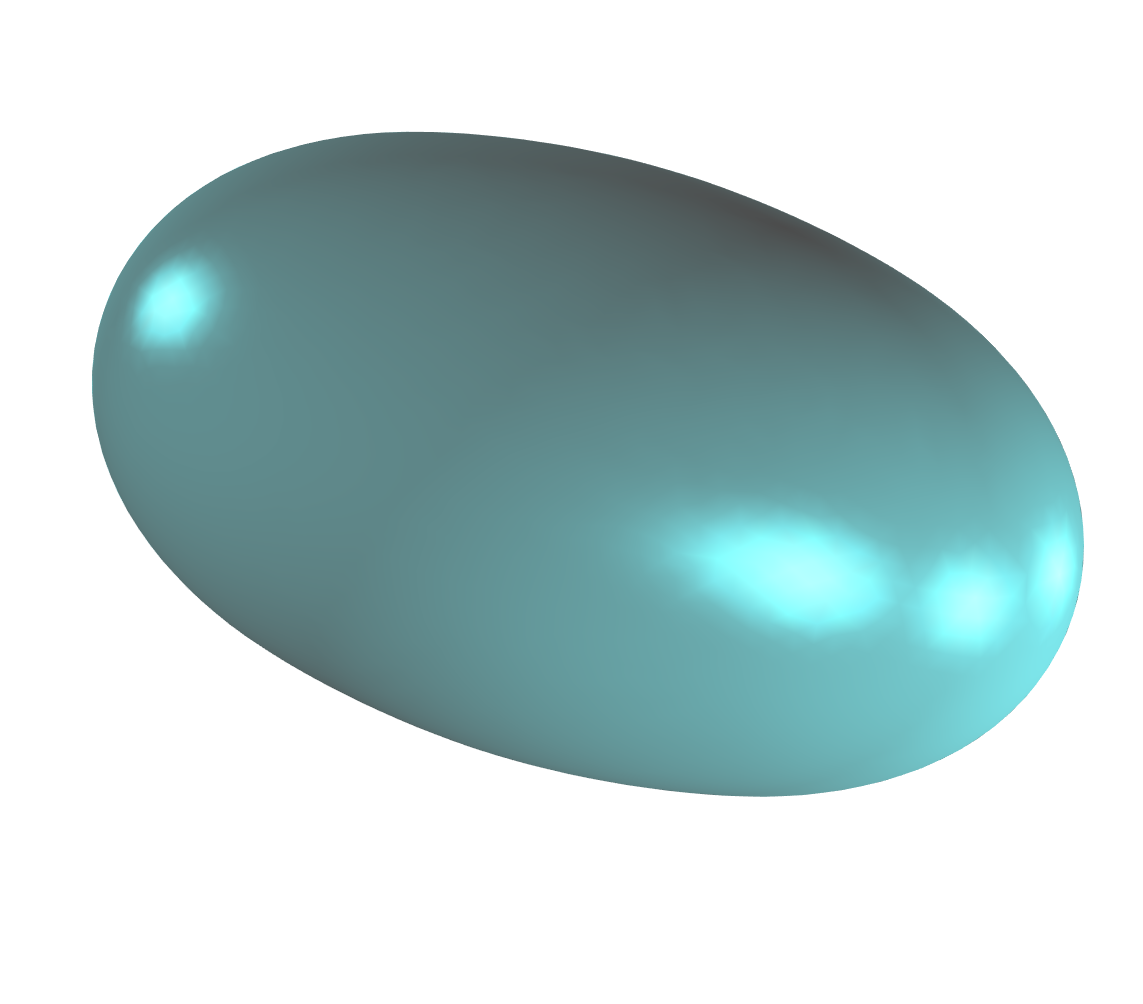}&
 \includegraphics[width=0.2\textwidth]{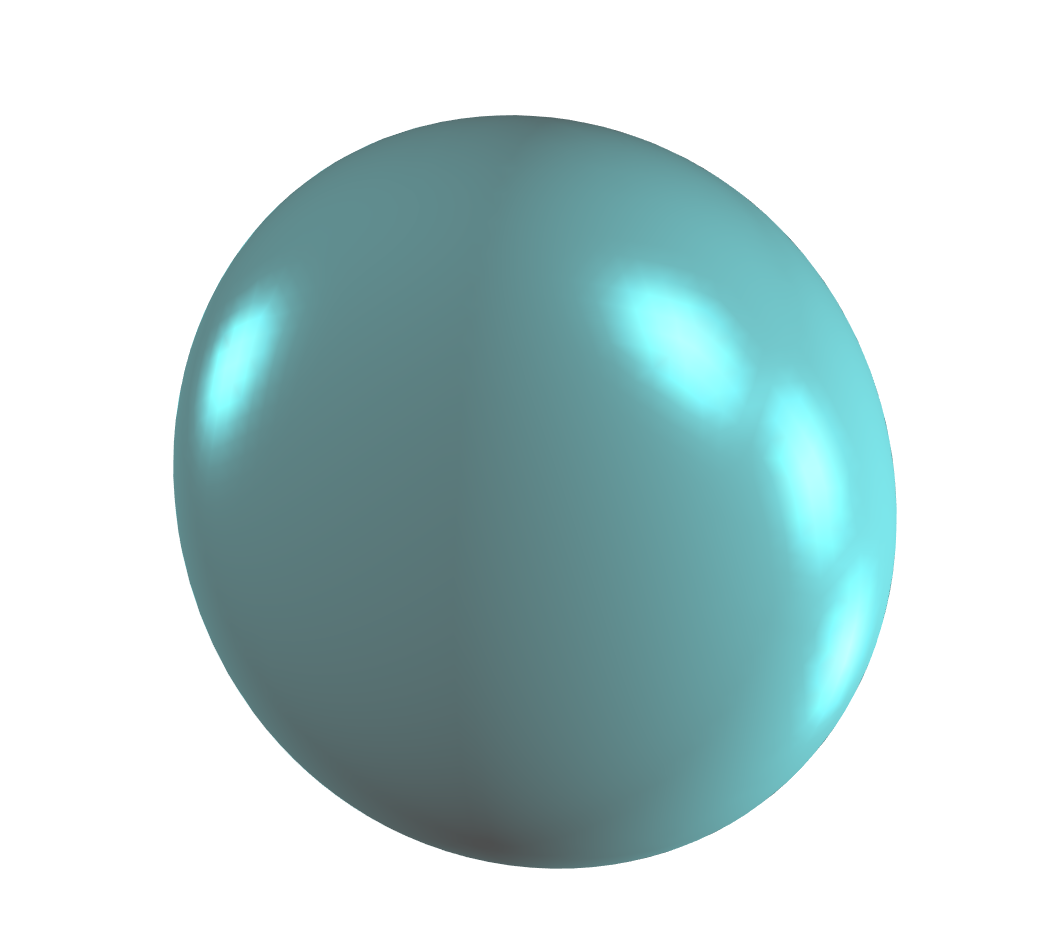} \\ 
 
  $\lambda_2=43.07$ & $\lambda_3=49.37$ & $\lambda_5=73.81$ & $\lambda_6=75.64$ \\
  
  \includegraphics[width=0.2\textwidth]{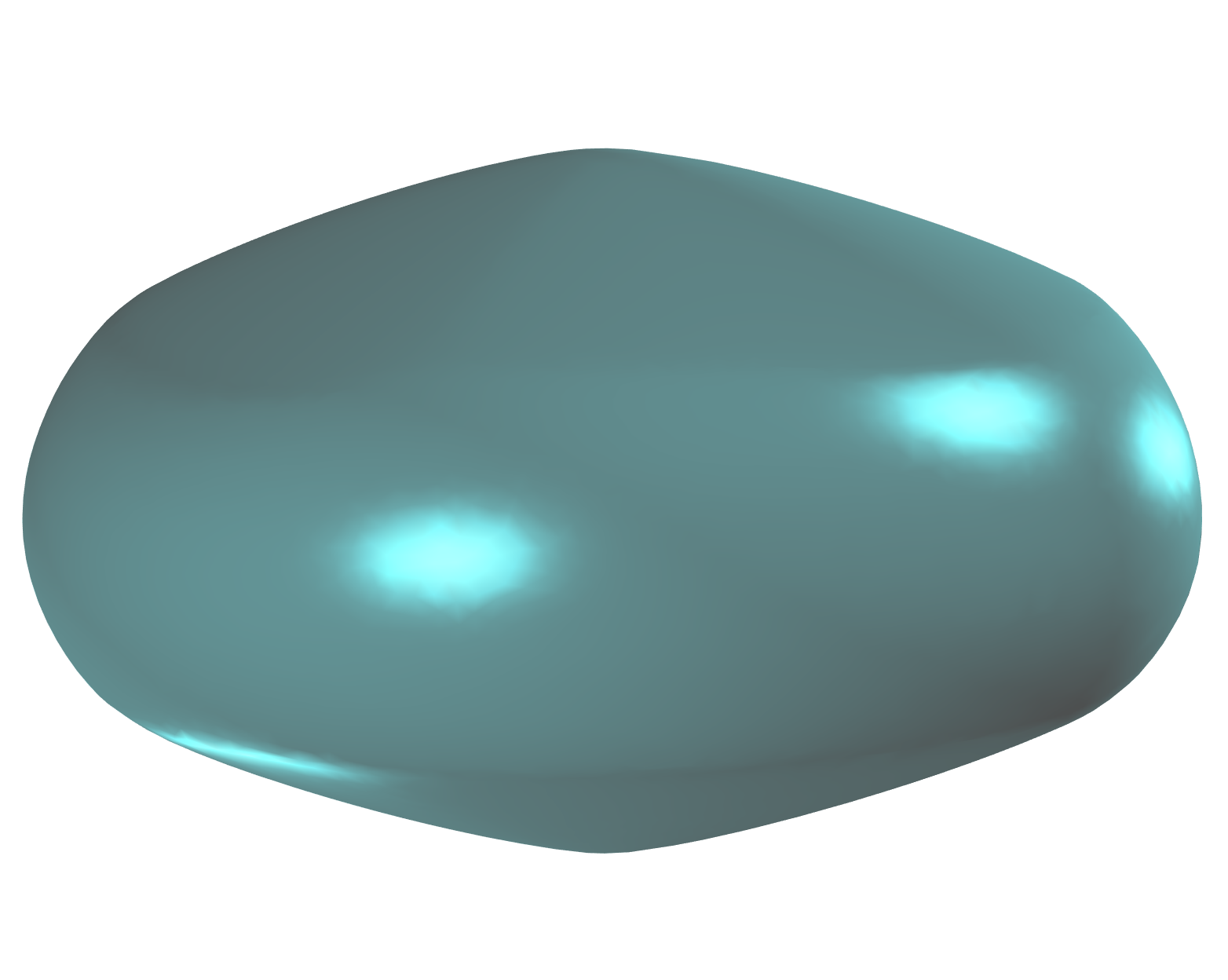}& 
  \includegraphics[width=0.2\textwidth]{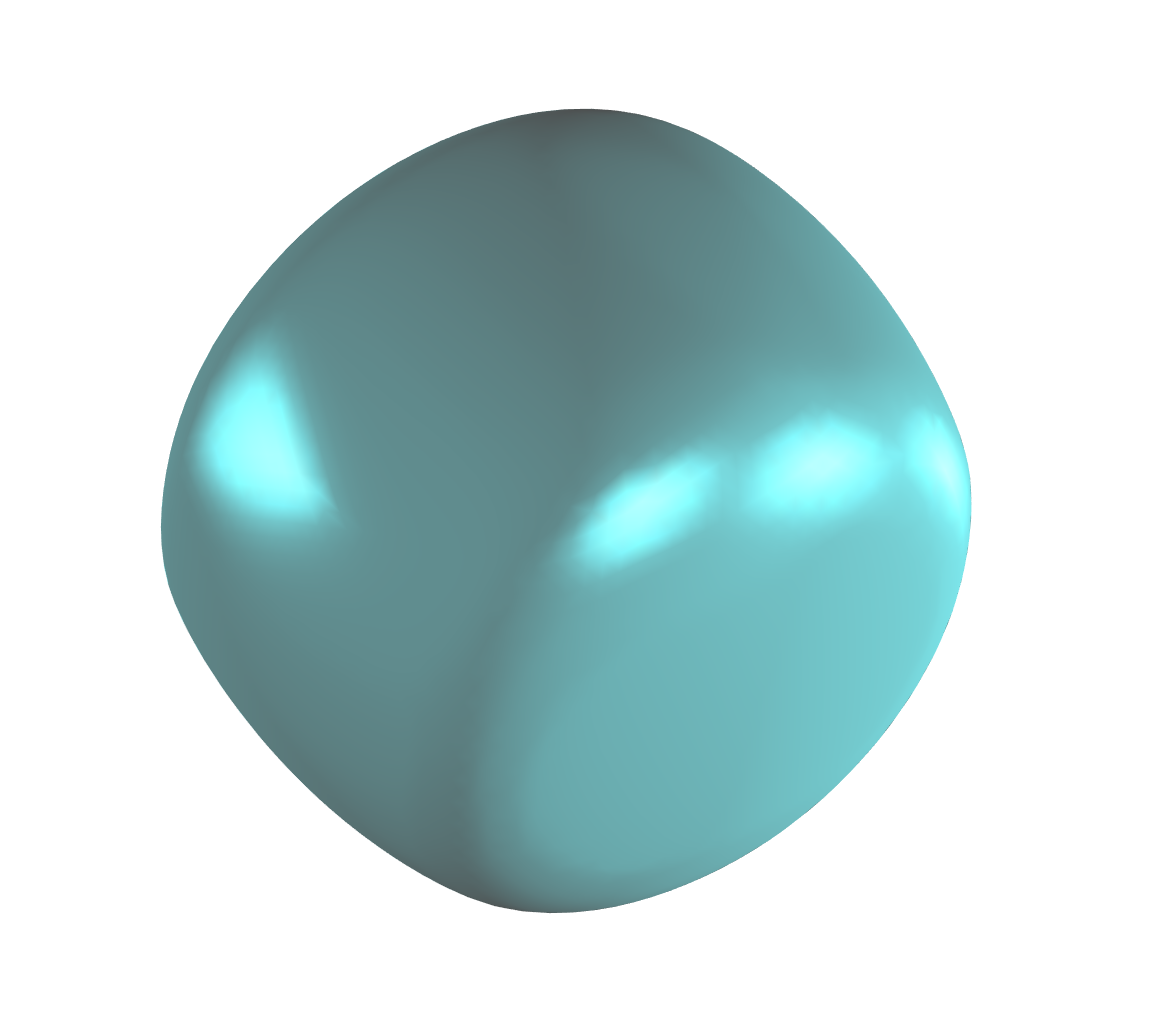}&
  \includegraphics[width=0.2\textwidth]{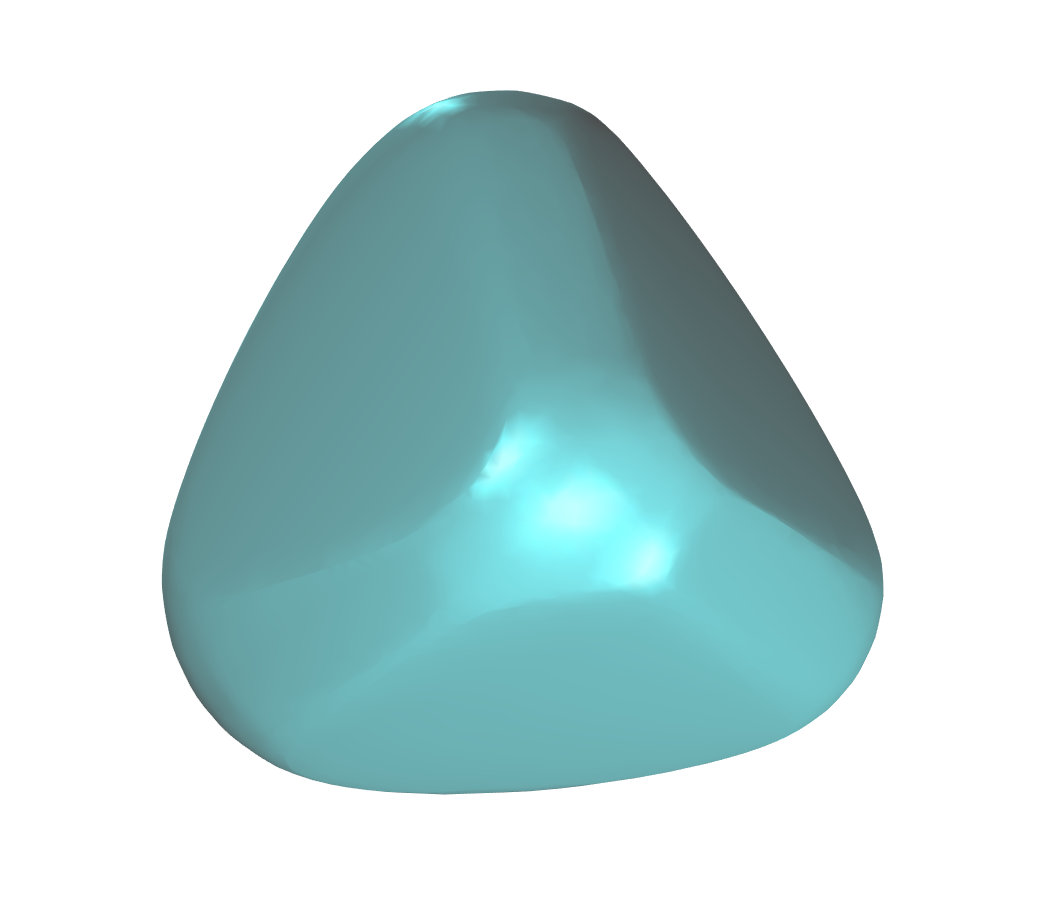} & \\
  $\lambda_7=80.74$ & $\lambda_8=85.29$ & $\lambda_{10}=93.22$ &
 \end{tabular}
 
 \caption{Optimization of eigenvalues under fixed volume and convexity constraint in 3D.  For $k \in \{4,9\}$ the ball is a numerical minimizer for $\lambda_k$ at fixed volume even without imposing the convexity constraint.}
 \label{fig:convexEig3D}
 \end{figure}

 The numerical discretization proposed in Section \ref{sec:dim3}  makes it possible to perform the same simulations in dimension three, with no additional difficulty. These results are shown in Figure \ref{fig:convexEig3D} for $k \in \{2,3,5,6,7,8,10\}$. For $k \in \{4,9\}$ balls are the numerical minimizers. It can be noticed that in some cases there are boundary regions of the numerical optimizers which seem to have at least one of the principal curvatures equal to zero. The numerical simulations are done using $900$ spherical harmonics and $M_c=5000$ points where convexity constraints are imposed. 
 
\subsection{Minimize the Dirichlet-Laplace eigenvalues under constant width constraint}
\label{sec:eigs-cw}

This section presents numerical approximations of solutions to Problem \ref{prob:eigs-cw} in dimension three. The two dimensional case was studied in \cite{BHL17}, where it was noted that the disk appears more often as a local minimizer. In fact, the precise list of indices $k$ for which the disk is a local minimizer for $\lambda_k(\omega)$ under constant width constraint in dimension two is given in \cite{BHL17}. 

\begin{figure}
	\centering
	\begin{tabular}{ccc}
		\includegraphics[height=0.27\textwidth]{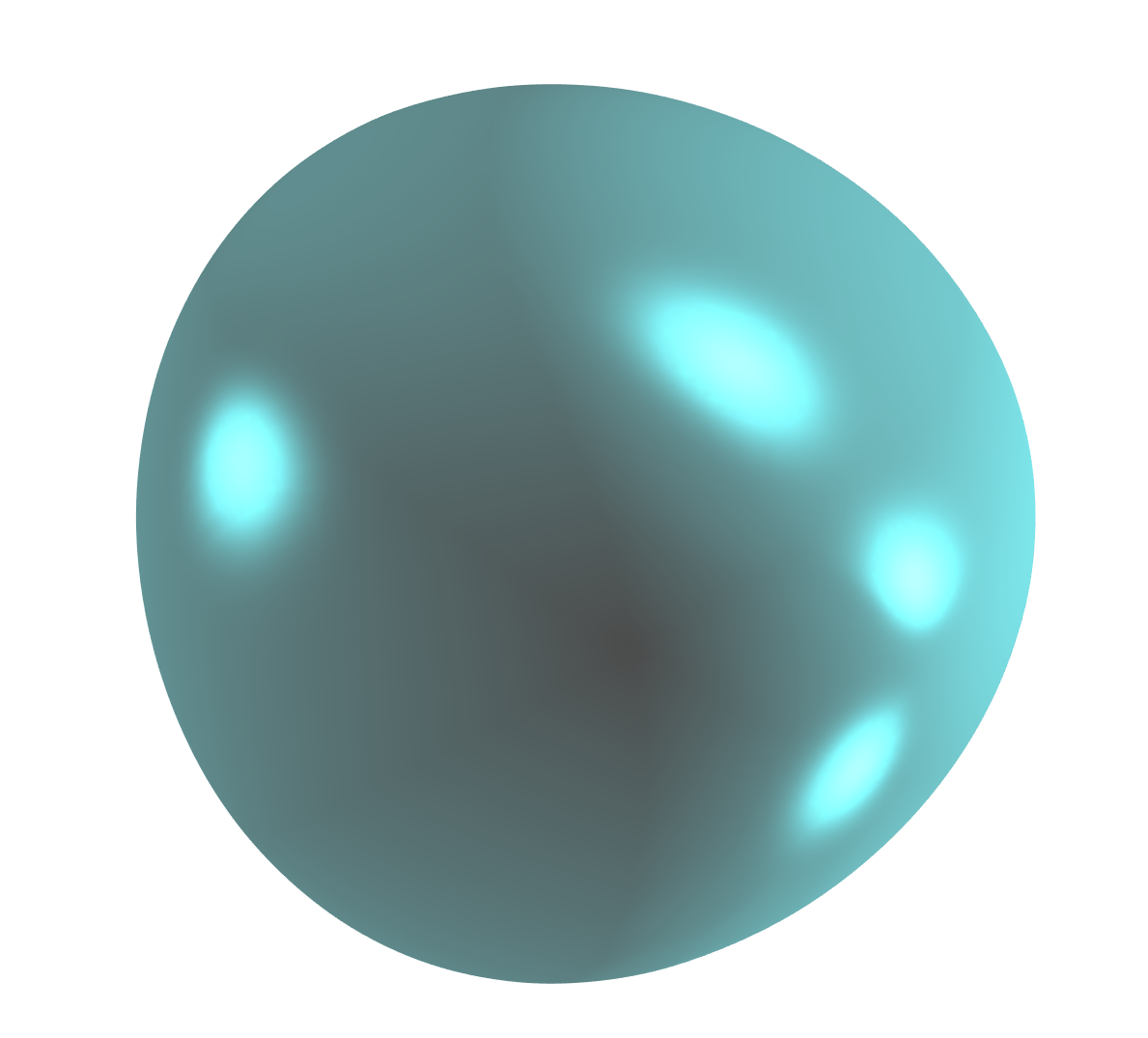}&
		\includegraphics[height=0.27\textwidth]{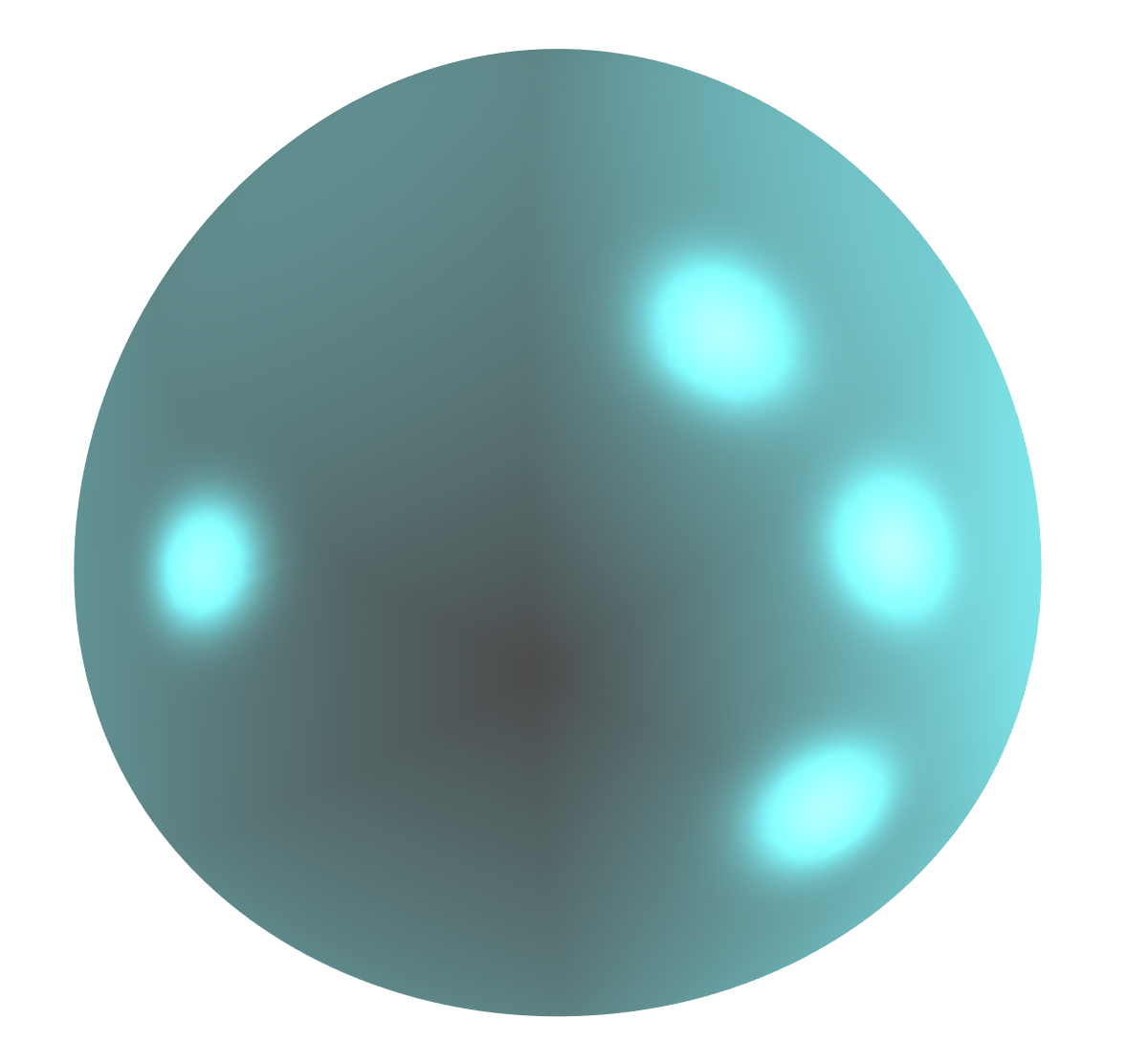}& 
		\includegraphics[height=0.27\textwidth]{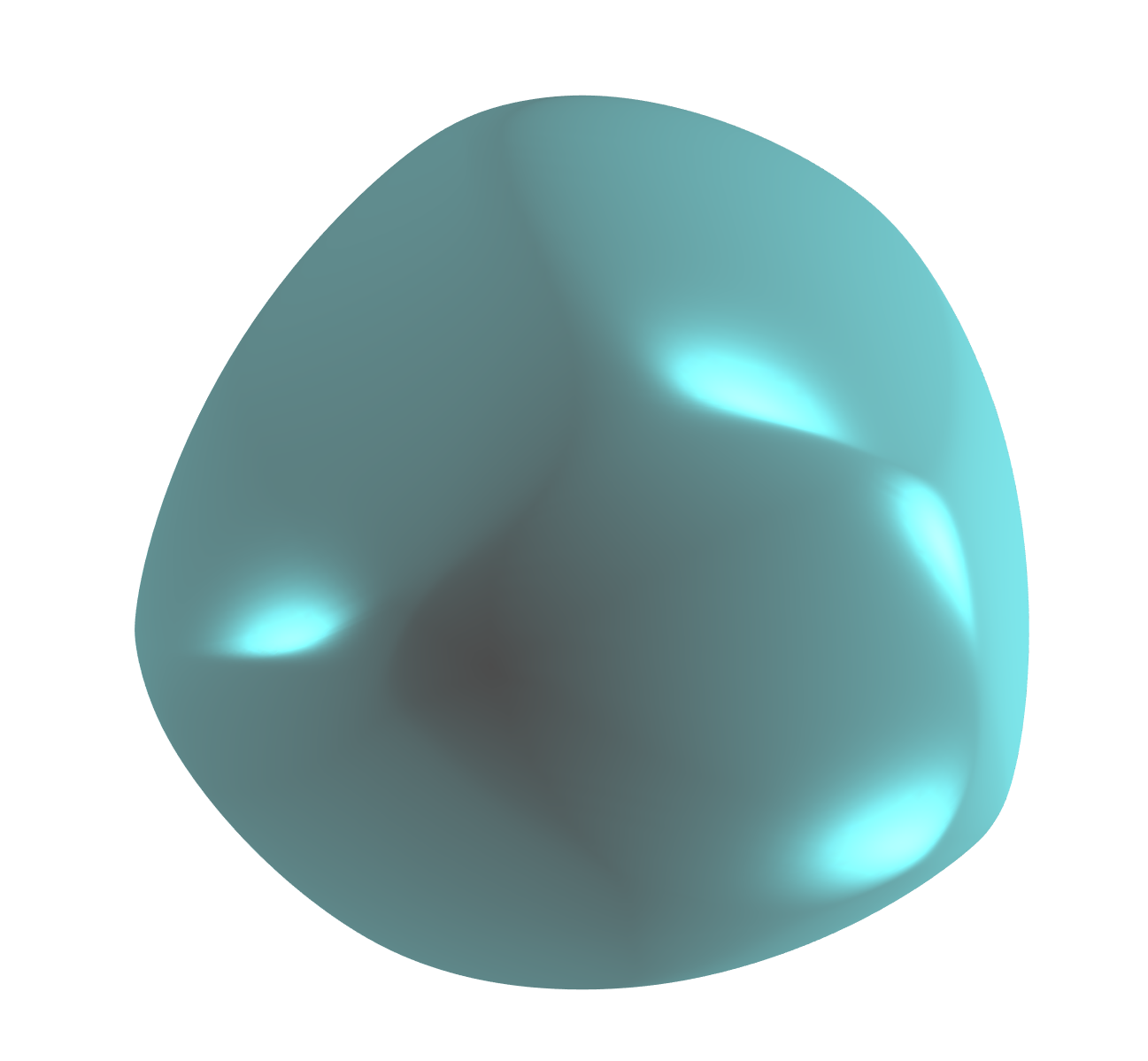} \\ 
		
		$\lambda_{10}=39.41$ & $\lambda_{46}=88.80$ &  $\lambda_{99}=156.14$ \\
		$\lambda_{10}(\bigcirc)=39.48$ & $\lambda_{46}(\bigcirc)=88.83$ &  $\lambda_{99}(\bigcirc)=157.92$ \\
	\end{tabular}
	\caption{Minimization of the eigenvalues $\lambda_k$, with $k\in\left\{10,46,99\right\}$ under fixed width in 3D.}
	\label{fig:cwEigs3D}
\end{figure}

The convexity and constant width constraints are imposed as indicated in Section \ref{sec:dim3}. We use gradient information in order to perform the optimization as indicated in Section \ref{sec:shderiv}. In our computations we observe that the ball appears often as a minimizer, but as observed in the two dimensional case in \cite{BHL17}, we expect that this only happens for finitely many indices $k$. Notable exceptions are the indices corresponding to a simple eigenvalue for the ball. Figure \ref{fig:cwEigs3D} shows the non-trivial shapes of constant width obtained with our algorithm for $k\in\left\{10,46,99\right\}$, the three smallest indices for which the corresponding eigenvalue of the ball is simple. For comparison the corresponding eigenvalue of the ball with the same width is shown in Figure~\ref{fig:cwEigs3D}. As in the two dimensional study in \cite{BHL17}, one could investigate the local minimality of the eigenvalues of the Dirichlet-Laplace operator on the ball in the class of constant width bodies. The numerical simulations use $900$ spherical harmonics and $3000$ points on where convexity constraints are imposed.

\subsection{Minimizing the volume under constant width constraint - the Meissner conjecture}
\label{sec:meissner}

This section presents numerical approximations of solutions to Problem \ref{prob:minvol-cw}. The convexity and constant width constraints are imposed as indicated in Section \ref{sec:numerics}. The area is explicit in terms of the Fourier coefficients in dimension two (see \eqref{eq:area-2d}). Moreover, the volume of constant width bodies is explicit in terms of the spherical harmonics coefficients in dimension three (see \eqref{eq:vol3D} and \eqref{eq:energy}).  Therefore, in this case the functional, its gradient and the corresponding Hessian matrix can be computed explicitly, leading to quickly converging numerical algorithms. 

\begin{figure}
\centering
\includegraphics[width=0.3\textwidth]{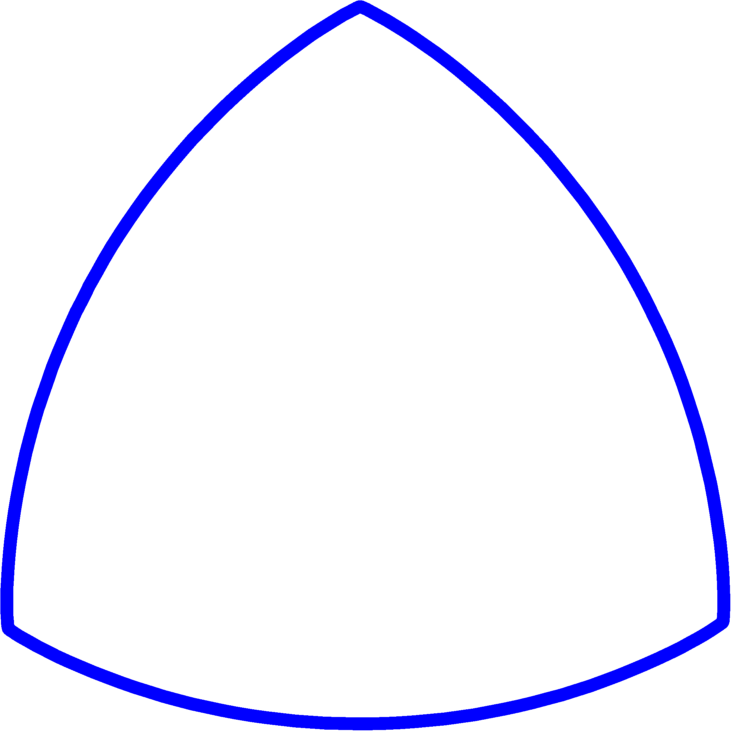}\quad \quad
\includegraphics[width=0.3\textwidth]{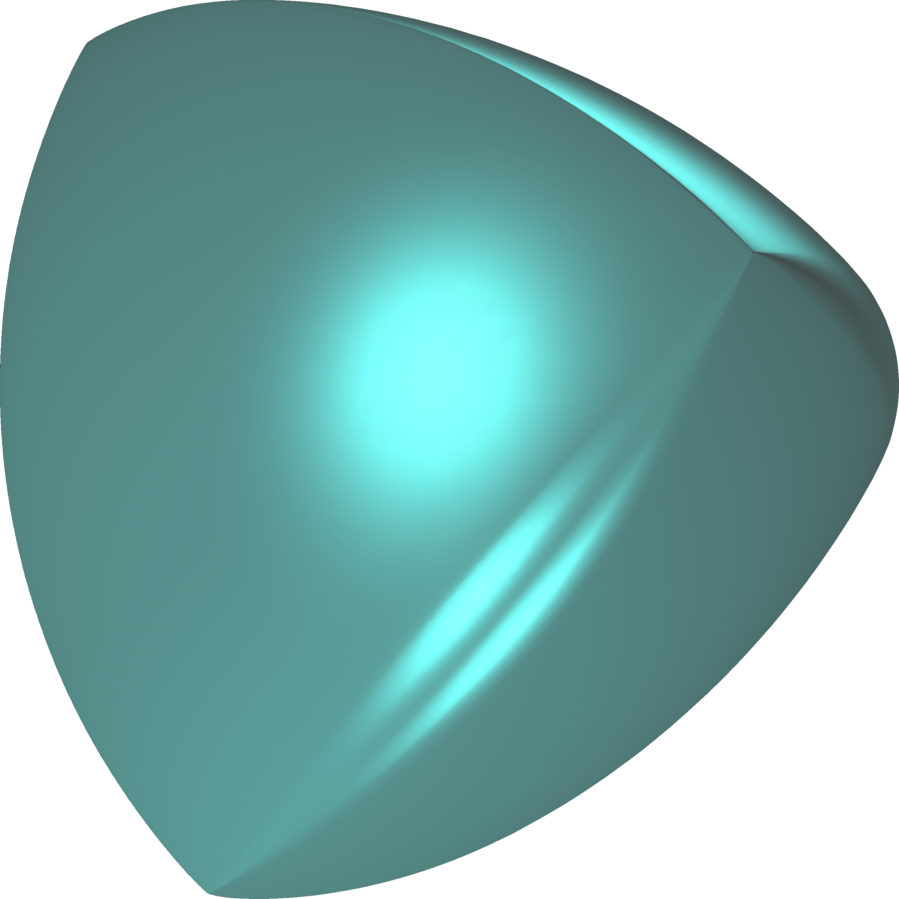}
\caption{Minimization of area and volume under constant width constraint: Reuleaux triangle (left) and Meissner's body (right)}
\label{fig:cwArea2D3D}
\end{figure}

The two dimensional result, the Reuleaux triangle, is shown in Figure \ref{fig:cwArea2D3D}. The minimal value for the area obtained with our algorithm, for $N=250$, corresponding to $501$ Fourier coefficients and width $w=2$ is $2.8196$. This is slightly larger than but very close to the explicit area of the Reuleaux triangle of width $2$ which is equal to $2(\pi-\sqrt{3}) = 2.8191$.

The minimization of the volume under constant width constraint in dimension three is a famous open problem. The conjectured optimizer is a Reuleaux tetrahedron with three rounded edges. There are two configurations, with the same volume, the difference being in the position of the rounded edges: all starting from one vertex or forming a triangle. These shapes are called Meissner's bodies. Various works deal with the analysis of 3D shapes of constant width which minimize the volume. Among these we cite \cite{KW11}, which presents many aspects related to the Meissner bodies and why they are conjectured to be optimal. It is mentioned that in \cite{muller_million} the author generates a million random three dimensional bodies of constant width, using techniques from \cite{spheroforms}. Among these many bodies of constant width, none had a smaller volume than the ones of Meissner. In \cite{oudetCW} the local optimality of the Meissner's body was verified using an optimization procedure with a different parametrization of constant width shapes. 

The approach we present below allows us to obtain the Meissner bodies as results of a direct optimization procedure, starting from random initializations. The formulas \eqref{eq:vol3D} and \eqref{eq:energy} allow us to write the volume as a quadratic expression of the coefficients of the spherical harmonics decomposition \eqref{sphDecomp}. As indicated by the Property \ref{prop:coef-cw}, the constant width condition is imposed by fixing the first coefficient and considering only odd spherical harmonics in the decomposition. The convexity condition is achieved by using a discrete version of \eqref{convexity3D}. We note that even though the convexity condition is non-linear, it is explicit enough such that we may compute its gradient. In this way, the minimization of the volume under constant width condition in dimension three becomes a constrained optimization problem of a quadratic functional with non-linear quadratic constraints. The optimization algorithm uses the explicit expression of the gradient and the Hessian matrix. For unit width, in view of \cite{KW11}, the Meissner bodies have volume equal to 
\[ V_M = \pi \left( \frac{2}{3}-\frac{\sqrt{3}}{4} \cdot \arccos\left( \frac{1}{3}\right) \right) = 0.419860.\]
In our computation, using $402$ spherical harmonics, i.e. using Legendre polynomials up to degree $28$, we obtain the shape represented in Figure \ref{fig:cwArea2D3D}, with volume $0.4224$. The shape obtained strongly resembles Meissner's body and its volume is about $0.6\%$ larger than $V_M$ presented above. This may be due to the fact that singularities in the surface of the Meissner bodies are not sufficiently well approximated using the number of spherical harmonics above. Starting from different random initial coefficients we always arrive at shapes which are close to one of the two Meissner bodies \cite{KW11}.

\subsection{Rotors of minimal volume}
\label{sec:rotors}
As underlined in Section \ref{sec:theory}, for some convex domains $P$ there exist convex shapes $\omega$ which can rotate inside $P$ while touching all its sides (or faces in dimension three). Therefore, when rotors exist, it makes sense to consider the problem of finding rotors of minimal volume. Theoretical aspects are recalled in the definition of Problem \ref{prob:rotors} and numerical computations for the two dimensional case are presented in \cite{BHcw12}. In dimension three there exist rotors for the cube (constant width bodies), the regular tetrahedron and the regular octahedron. A characterization of rotors in terms of the coefficients of the decomposition of the support function can be found in \cite{goldberg_rotors}. More precisely we have the following:
\begin{itemize}
	\item In dimension two every regular polygon admits non-circular rotors. If the regular polygon has $n$ sides, $n \geq 3$, then only the Fourier coefficients of the support function for which the index has the form $nq\pm 1$ are non-zero, where $q$ is a positive integer.
	\item The rotors in a cube are bodies of constant width.
	\item The support function of a rotor in a regular tetrahedron has non-zero coefficients for the spherical harmonics of indices $0, 1,2$ and $5$
	\item The support function of a rotor in a regular octahedron has non-zero coefficients for the spherical harmonics of indices $0,1$ and $5$.
\end{itemize}  
The constant term in the spectral decomposition of the support function of a rotor in $P$ corresponds to the inradius of $P$. We note that when taking the midpoints of the edges of the regular tetrahedron we obtain a regular octahedron with the same inradius. Therefore rotors in the octahedron are also rotors for the tetrahedron which was already apparent from the characterization using the spherical harmonic coefficients.

\begin{figure}
\centering
\includegraphics[width=0.2\textwidth]{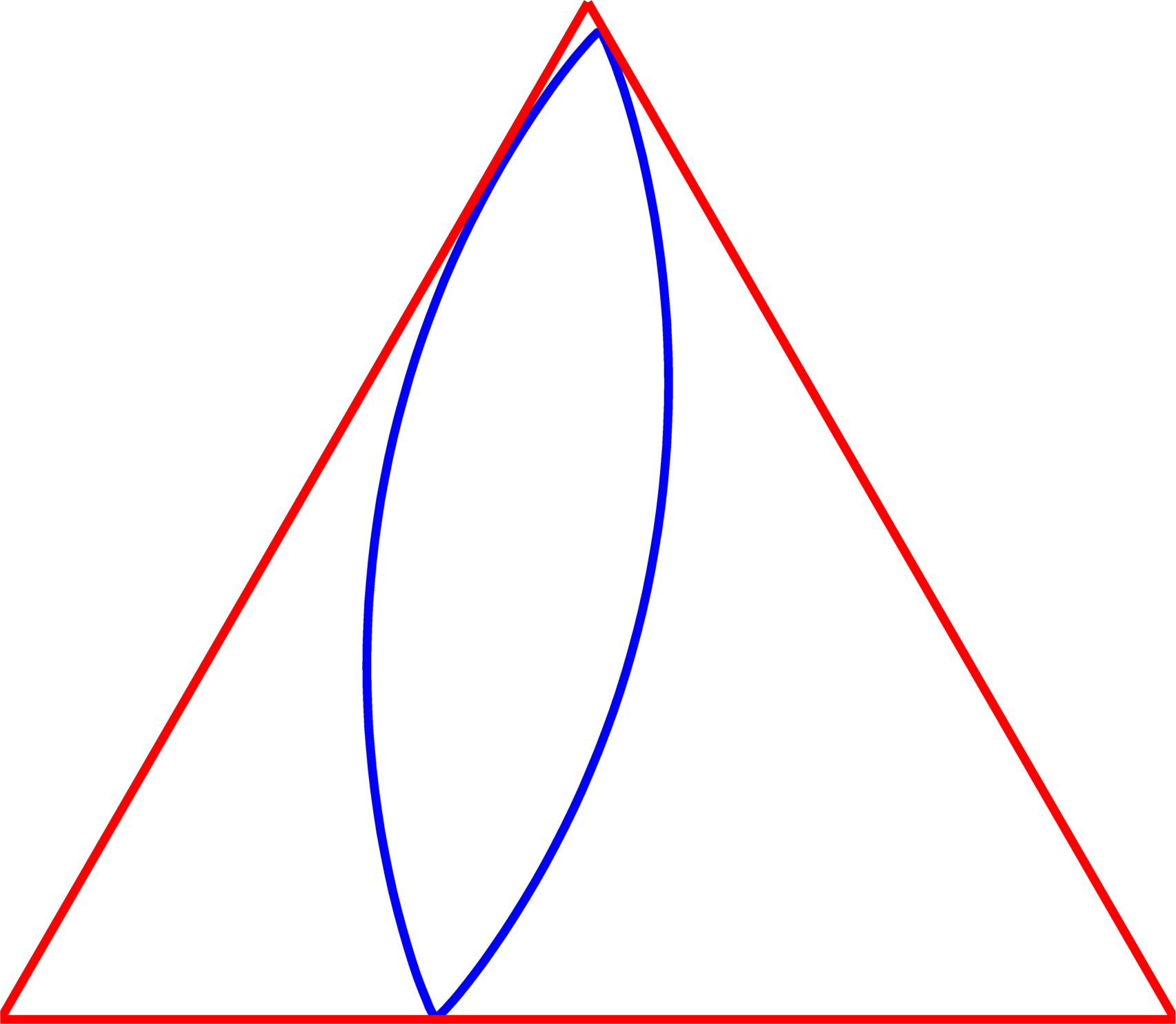}\quad 
\includegraphics[width=0.2\textwidth]{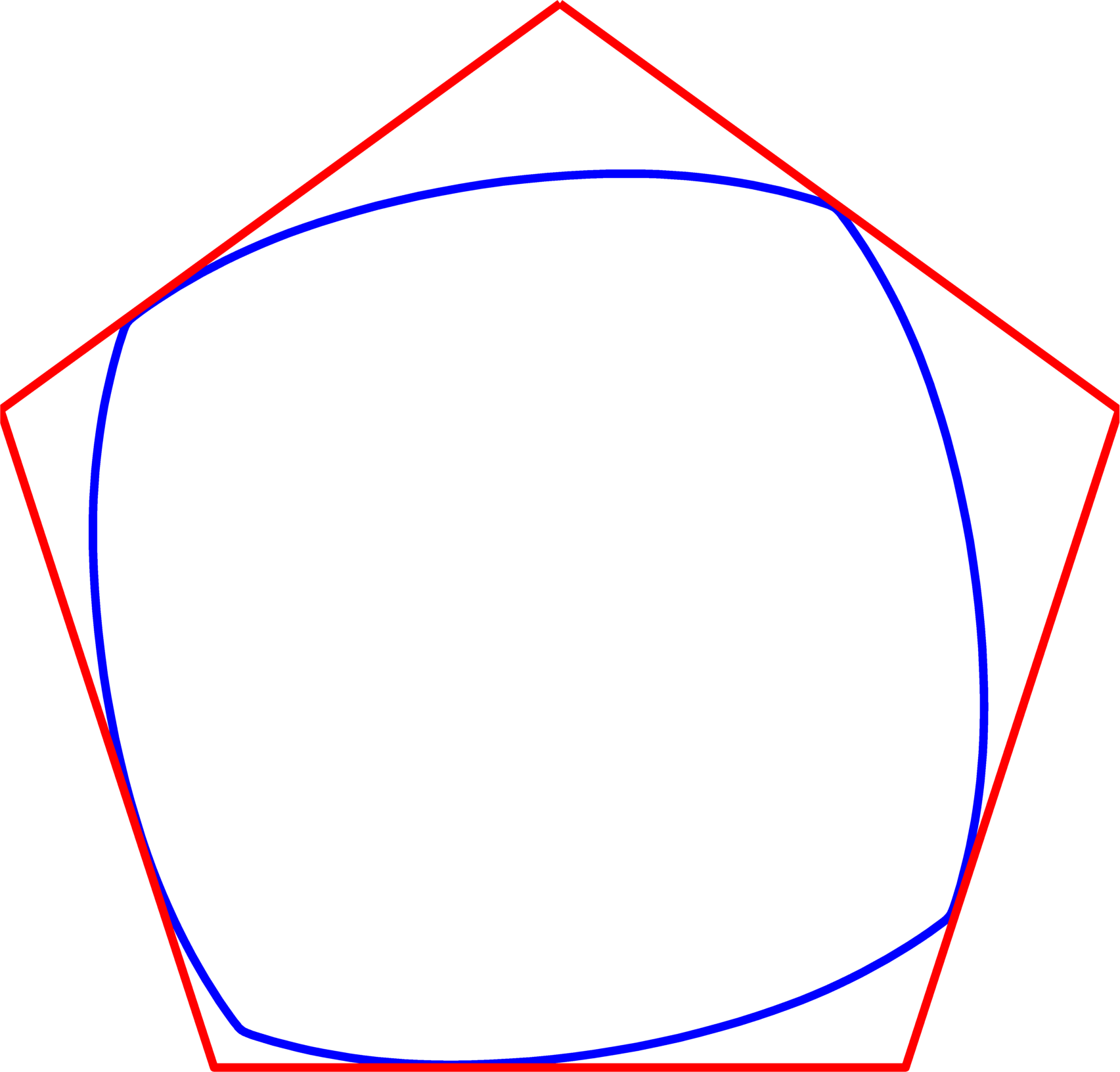}\quad 
\includegraphics[width=0.2\textwidth]{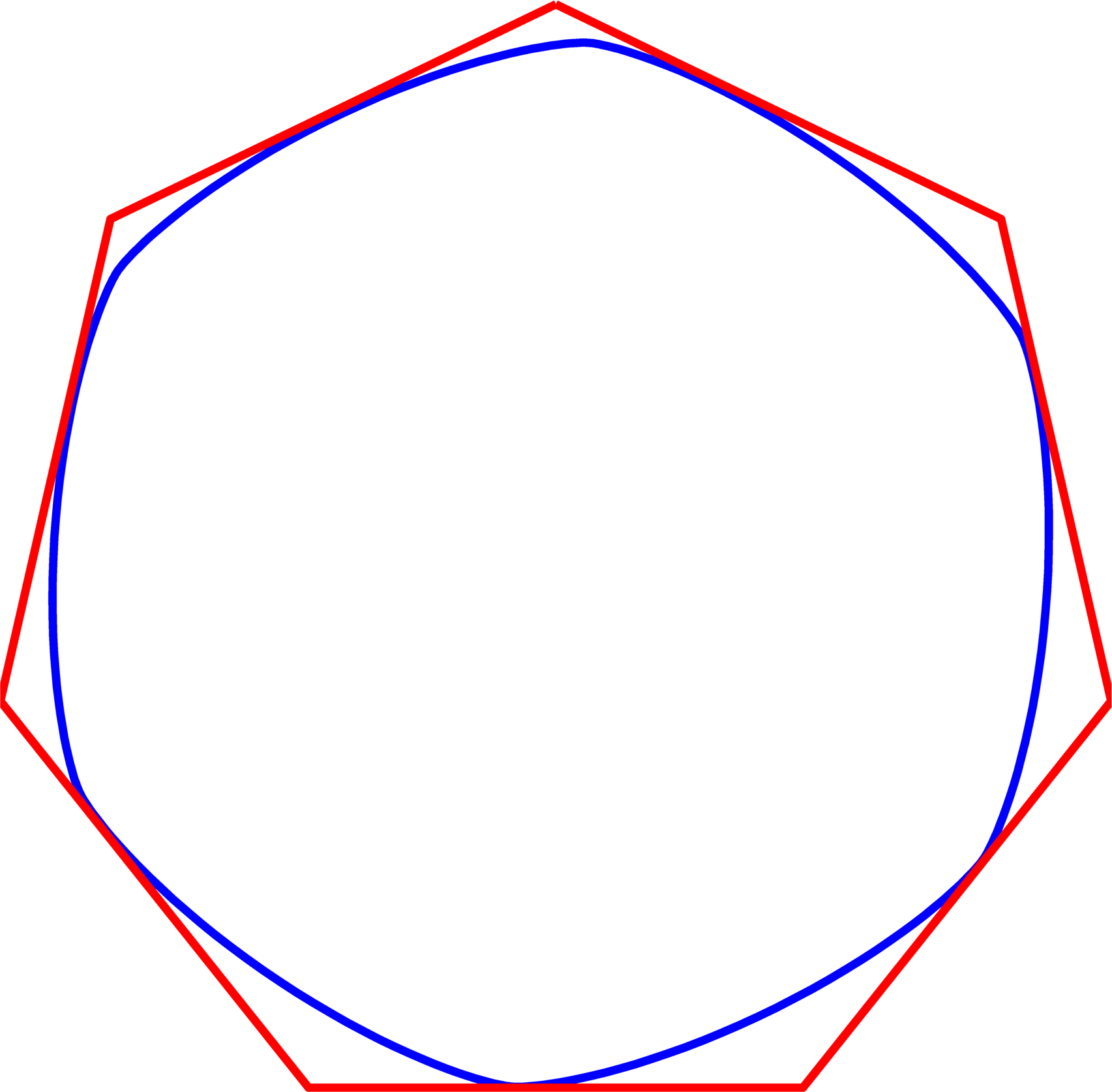} 
\caption{Examples of minimal area rotors in dimension two.}
\label{fig:rotors2D}
\end{figure}

\begin{figure}
\centering
\includegraphics[width=0.2\textwidth]{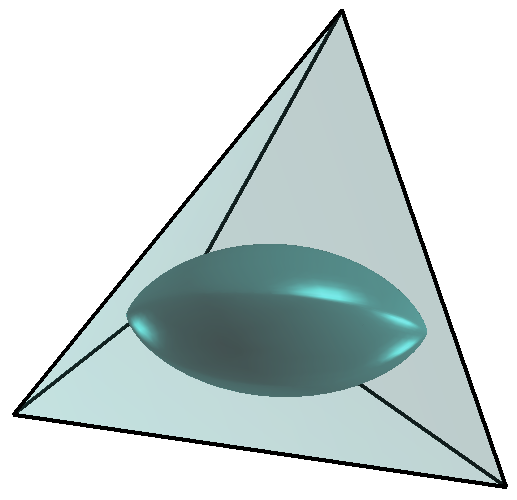}\quad
\includegraphics[width=0.18\textwidth]{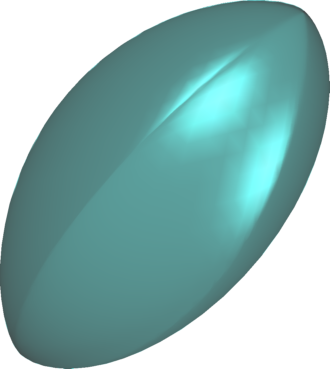}\quad 
\includegraphics[width=0.2\textwidth]{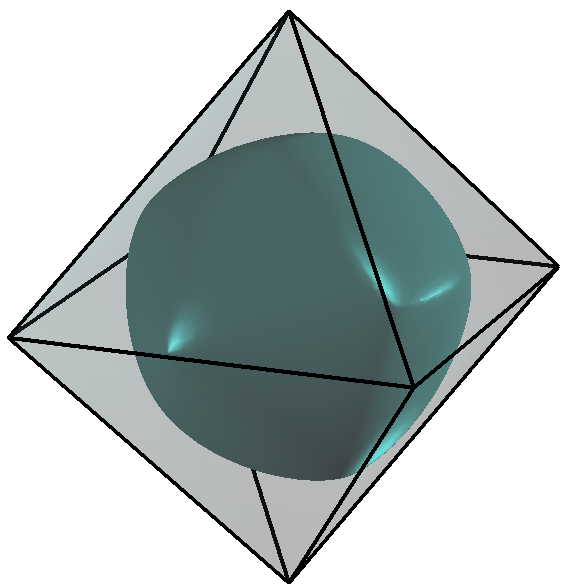}\quad
\includegraphics[width=0.2\textwidth]{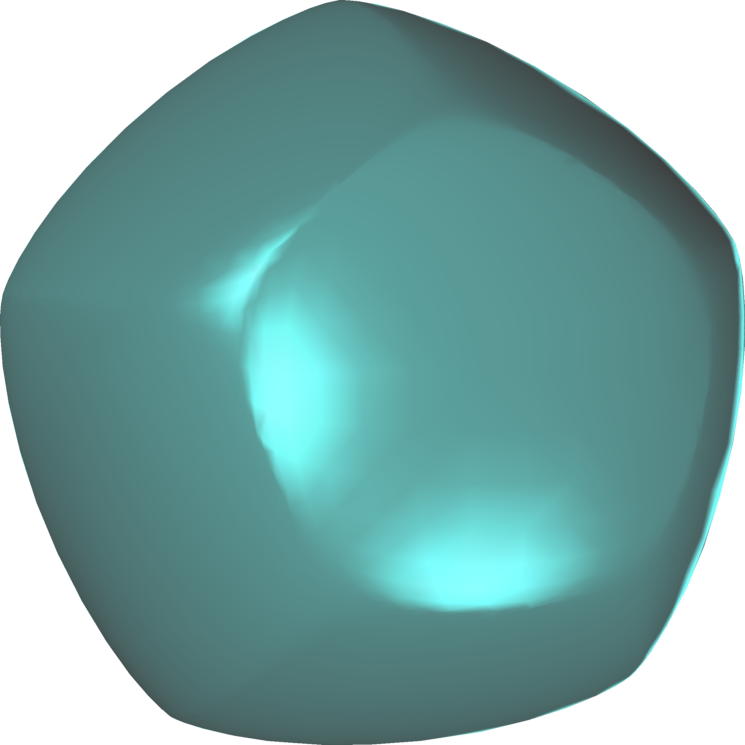}\quad 
\caption{Minimal volume rotors in the regular tetrahedron and the regular octahedron. Volume of the inscribed sphere: $0.5236$. Volume of the rotors: tetrahedron $0.3936$, octahedron $0.5041$.}
\label{fig:rotors3D}
\end{figure}

Computations of optimal rotors in dimension two were also made in \cite{BHcw12}, while the computations in dimension three are new. Note that rotors of maximal area and volume are the inscribed disc and the inscribed ball, respectively. Some two dimensional rotors of minimal area are shown in Figure \ref{fig:rotors2D}. Minimal volume rotors obtained numerically for the regular tetrahedron and the regular octahedron are shown in Figure \ref{fig:rotors3D}. In each case we consider an optimization problem depending only on the non-zero coefficients describing the rotors and we impose discrete convexity constraints like in Section \ref{sec:numerics}. 

The computations presented in Figure \ref{fig:rotors3D} are made for solids with inradius equal to $0.5$, corresponding to an inscribed ball of diameter $1$. Compared to the volume of the ball $B$ with unit diameter which is equal to $\pi/6=0.5236$ the minimal volume found numerically of a rotor in the tetrahedron and the octahedron circumscribed to the same ball $B$ are $0.3936$ and $0.5041$, respectively. Numerical minimizers for the tetrahedron and octahedron seem to be symmetric under a rotation of angle $2\pi/5$. This is due to the fact that the only coefficients which may change the geometry of the rotors correspond to the spherical harmonics of order $2$ or of order $5$. We recall that changing coefficients for spherical harmonics of order $1$ corresponds to translations. In particular, when searching for minimal rotors in the tetrahedron using only spherical harmonics of order $1$ and $2$ we get a radially symmetric minimizer of volume $0.4024$ which is slightly larger than the result including spherical harmonics of order $5$.

\subsection{Minimal width constraint}

Problem \ref{prob:area-minw} was first considered in \cite{oudetCW}. The area of a three-dimensional convex body is minimized under a minimal width constraint. The minimal width constraint is characterized by the inequality $p(\theta)+p(-\theta) \geq 1$, for every $\theta \in \Bbb{S}^2$.

\begin{figure}
\centering
\includegraphics[width=0.3\textwidth]{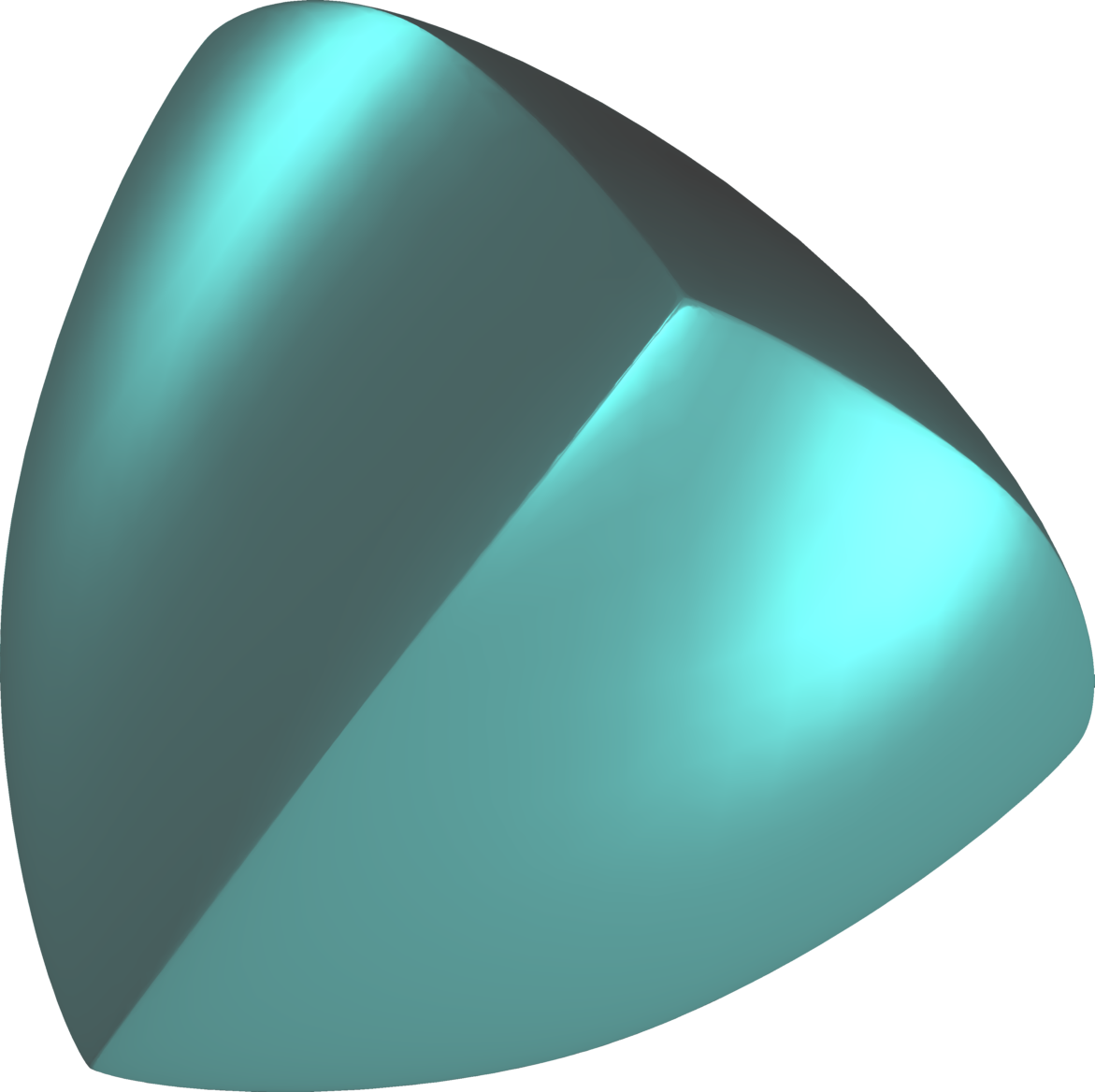}
\caption{Optimization under diameter bounds in dimension three. Minimization of the area for shapes having width at least $1$. The minimal area found by our algorithm is $2.9154$. }
\label{fig:diam3D}
\end{figure}

 As underlined in Section \ref{sec:dim3}, diameter bounds can be imposed in an approximate way using a finite number of points uniformly distributed on the sphere. This is done in the same way as the discrete convexity condition. These bounds on the diameter correspond to a set of linear inequality constraints. In Figure \ref{fig:diam3D} we present the result given by the algorithm. The shape resembles the optimizer given in \cite{oudetCW} and the value of the functional is slightly improved. In the computations $M_c=2000$ points are used for the discrete convexity condition and $M_d=1000$ pairs of opposite diametral points for computing the discrete diameter inequalities. We used $250$ spherical harmonics in the decomposition of the support function. The computation of the area is explicit in terms of the spherical harmonics coefficients, as shown in \eqref{eq:area3D}. The minimal area obtained with our algorithm is $2.9154$ which is slightly smaller than $2.9249$, the value of the minimal area in the result presented in \cite{oudetCW}.

\subsection{Inclusion constraint}

In this Section we show how to impose inclusion constraints for shape optimization problems. As recalled in Section \ref{support}, two convex bodies $B_1,B_2$ in $\Bbb{R}^n$, with support functions $p_{B_1},p_{B_2}$, respectively satisfy the inclusion constraint $B_1 \subset B_2$ if and only if the support functions verify 
\begin{equation}
p_{B_1}(\theta) \leq p_{B_2}(\theta) \text{  for every  }\theta \in \Bbb{S}^{n-1}.
\label{eq:incl_ineq}
\end{equation}
 As in the case of the convexity and diameter constraints, we impose \eqref{eq:incl_ineq} on a sufficiently dense discrete subset of $\Bbb{S}^{n-1}$. In dimension three, when dealing with Cheeger sets for polyhedra, it is enough to impose the inclusion constraints only for directions which are normals to the faces of the polyhedron. This simplifies the optimization algorithm by decreasing the number of constraints.
 
\begin{figure}
\centering
\includegraphics[height=0.22\textwidth]{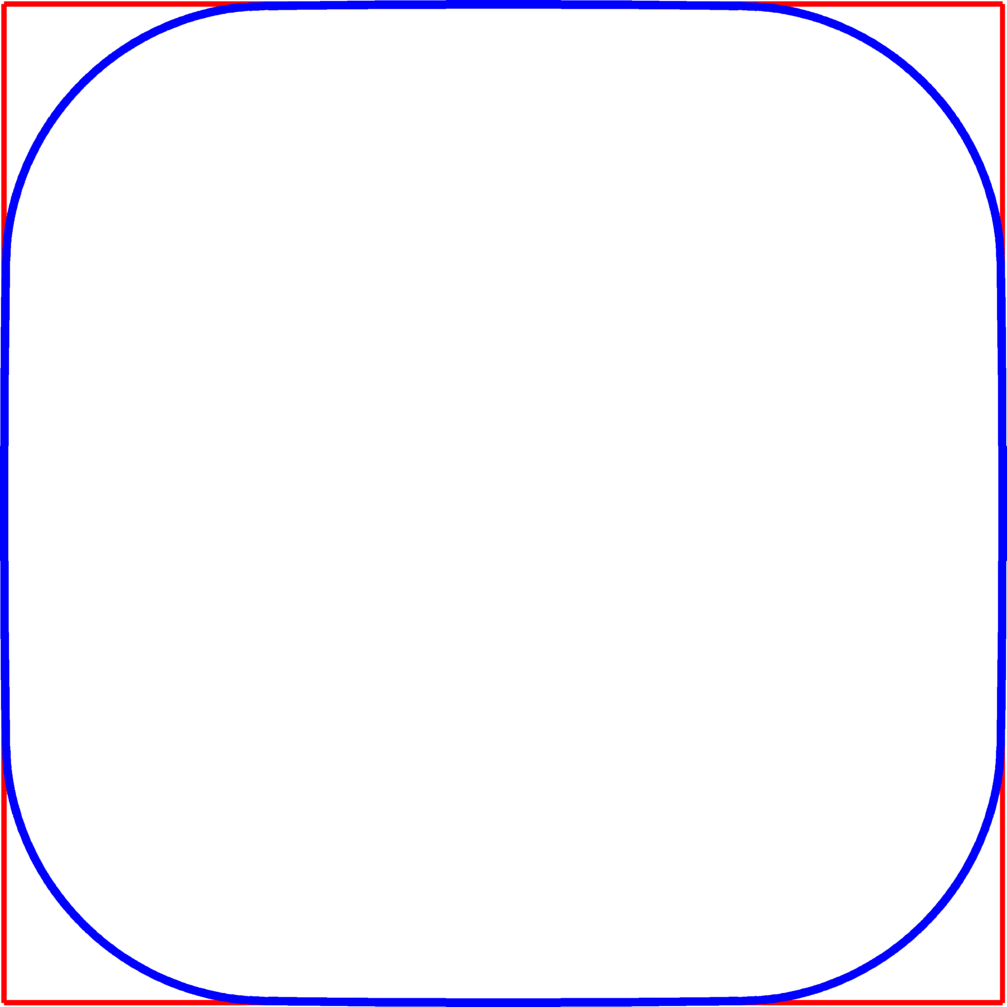}\quad
\includegraphics[height=0.22\textwidth]{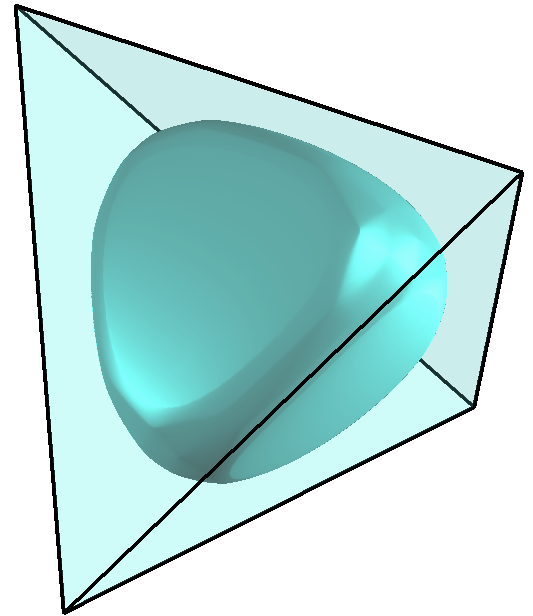}\quad
\includegraphics[height=0.22\textwidth]{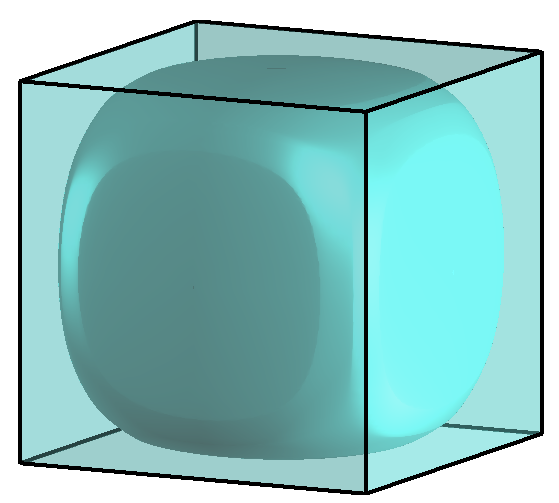}\quad
\includegraphics[height=0.22\textwidth]{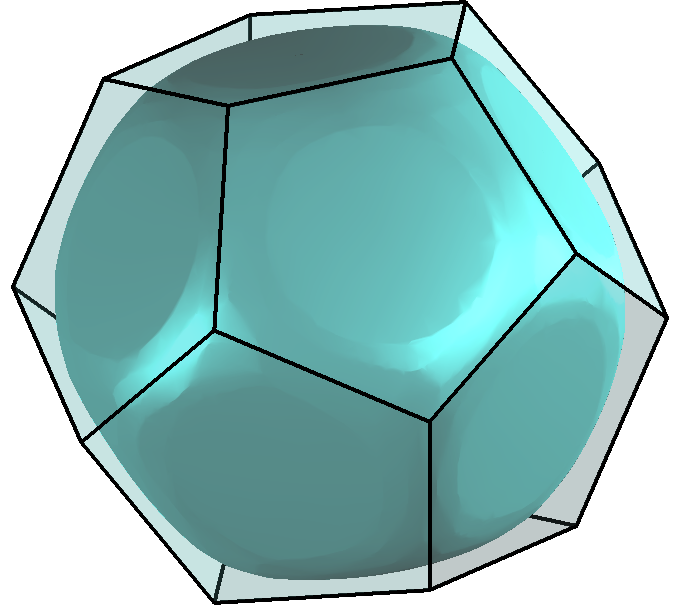}
\caption{Computation of Cheeger sets by optimizing the ratio perimeter/volume under convexity and inclusion constraints.}
\label{fig:Cheeger}
\end{figure}

 Cheeger sets provide a classical example of shape optimization problems involving convexity and inclusion constraints. The theoretical formulation is given in Problem \ref{prob:Cheeger}. The Cheeger sets of $\Omega$ minimize the ratio perimeter/area (or surface area/volume) for convex subsets of $\Omega$. This objective function can be computed and optimized using the proposed algorithm. As shown in Section \ref{sec:numerics}, convexity and inclusion constraints are discretized as linear inequalities in terms of the coefficients of the spectral decomposition of the support function. Some examples of computation of Cheeger sets for the square in the plane and for the regular tetrahedron, the cube and the regular dodecahedron in dimension three are shown in Figure \ref{fig:Cheeger}. 

\section{Conclusions}

In this work, the properties of the support function are used to deal numerically with various constraints in shape optimization problems. The spectral decomposition of the support function using Fourier series in dimension two and spherical harmonics in dimension three are particularly well suited in order to discretize convexity, constant-width, diameter and inclusion constraints. The numerical tests use standard tools readily available in optimization software like quasi-Newton or Newton methods with linear or non-linear constraints and cover a wide variety of shape optimization problems with various constraints.

\bibliography{./supportf_applications}
\bibliographystyle{./sn-mathphys}

\end{document}